\documentclass[12pt]{amsart}
\usepackage[centertags]{amsmath}
\usepackage{amsfonts}
\usepackage{amssymb}
\usepackage{amsthm}
\usepackage{newlfont}
\DeclareMathOperator{\Int}{Int}
\DeclareMathOperator{\HPath}{HPath}
\DeclareMathOperator{\HSkel}{HSkel}
\DeclareMathOperator{\Co}{Co}
\DeclareMathOperator{\dist}{dist}
\DeclareMathOperator{\diam}{diam}
\DeclareMathOperator{\MLR}{MLR}

\begin{document}

\theoremstyle{plain}
\newtheorem{thm}{Theorem}[section]
\newtheorem{cor}[thm]{Corollary}
\newtheorem{lem}[thm]{Lemma}
\newtheorem{prop}[thm]{Proposition}
\newtheorem{Theorem}{Theorem}[section]
\newtheorem{Corollary}[Theorem]{Corollary}
\newtheorem{Lemma}[Theorem]{Lemma}
\theoremstyle{definition}
\newtheorem{defn}{Definition}[section]
\newtheorem{Definition}[Theorem]{Definition}
\newtheorem{Remark}[Theorem]{Remark}
\newtheorem{Main Theorem}{Main Theorem}
\theoremstyle{remark}
\newtheorem{rem}{Remark}[section]

\numberwithin{equation}{section}
\renewcommand{\theequation}{\thesection.\arabic{equation}}

\title[Lower Bounds for Boundary Roughness]
  {Lower Bounds for Boundary Roughness\\
  for Droplets in Bernoulli Percolation}
\author{Hasan B. Uzun}
\address{ALEKS Corporation\\
  400 N. Tustin Avenue, Suite 300\\
  Santa Ana, CA 92705  USA}
\email{huzun@aleks.com}
\author{Kenneth S. Alexander}
\address{Department of Mathematics DRB 155\\
  University of Southern California\\
  Los Angeles, CA  90089-1113 USA}
\email{alexandr@math.usc.edu}
\thanks{The research of the first author was supported by NSF grant
DMS-9802368. The research of the second author was supported by NSF
grants DMS-9802368 and DMS-0103790.}

\date{\today}

\begin{abstract}
We consider boundary roughness for the ``droplet'' created when
supercritical two-dimensional Bernoulli percolation is conditioned to
have an open dual circuit surrounding the origin and enclosing an area
at least $l^2$, for large $l$.  The maximum local roughness is the
maximum inward deviation of the droplet boundary from the boundary of
its 	own convex hull; we show that for large $l$ this maximum is at
least of order $l^{1/3}(\log l)^{-2/3}$.  This complements the upper
bound of order $l^{1/3}(\log l)^{2/3}$ proved in \cite{Al3} for the
average local roughness.  The exponent 1/3 on $l$ here is in keeping with
predictions from the physics literature for interfaces in two
dimensions.
\end{abstract}
\maketitle

\section{Introduction}
We consider Bernoulli bond percolation on the square lattice at
supercritical density, conditioned to have a large dual
circuit enclosing the origin; we denote the outermost such circuit by
$\Gamma_0$.  (Complete definitions  and the basic properties of the
model will be given in the next section.)  The supercritical, or
percolating, regime of Bernoulli percolation is the analog of the
low-temperature phase of a spin system, and the region enclosed by
the dual circuit is the analog of the droplet that occurs with high
probability in the Ising magnet below the critical temperature
in a finite box with minus boundary condition, when it is
conditioned to have a number of plus spins somewhat larger than is
typical \cite{DKS}. In fact, the droplet boundary in the Ising
magnet appears as a circuit of open dual bonds in the corresponding
Fortuin-Kastelyn random cluster model (briefly, the FK model) of
\cite{FK}, in view of the construction given in
\cite{ES}.  One can gain information for
the study of the Ising droplet by studying the FK model conditioned
on $\Gamma_{0}$ enclosing at least a given area $l^{2}$, as is
done in \cite{Al3}. The droplet boundary in this FK model thus
corresponds to an interface; the heuristics in the case of Bernoulli
percolation are the same, but the mathematics is more tractable.  
We therefore refer to $\Gamma_0$ and its interior as a droplet.
Our main result is a lower bound on the
maximum local roughness of the droplet, that is, the maximum inward
deviation of the boundary of the droplet from the boundary of its
convex hull.  Related upper bounds were proved in \cite{Al3}.

\par The study of the shapes of such droplets is related to a
classical problem: When a fixed volume of one phase is
immersed in another, what is the equilibrium shape of the droplet, or
crystal, having minimal surface tension?  When the surface tension is
known, this is an isoperimetric problem. The solution of the
continuum version of the problem is given by Wulff \cite{Wu}: Let
$\tau(\textbf{n})$ be the surface tension of a flat interface
orthogonal to the outward normal $\textbf{n}$.  For a fixed crystal
volume, the equilibrium shape is given by the convex set
\begin{eqnarray} \label{E: classicalW}
W= \{ \textbf{x} \in \mathbf{R}^{d} \ | \ \textbf{x} \cdot
\textbf{n} \leq \tau(\textbf{n}), \ \mbox{for all } \textbf{n} \}.
\end{eqnarray}
In the two-dimensional Ising model, say with minus boundary
condition and conditioned to have an excess of pluses, a rigorous
justification of the Wulff construction has been given for the
resulting droplet of plus phase.  Minlos and Sinai considered an
instance in which the temperature $T$ tends to zero as the volume
grows to infinity, and proved that most of
the excess plus spins form a single droplet of essentially square
shape (\cite{MS1}, \cite{MS2}); the Wulff shape $W$ also tends to a
square as $T
\to 0$. Dobrushin, Kotecky and Shlosman \cite{DKS} then provided a
justification of the  Wulff construction at very low
fixed temperatures.  Moreover, they showed that the
Hausdorff distance between the droplet boundary $\gamma$ and the
boundary of the Wulff shape $W$ is bounded by a power of the
linear scale of the droplet.  This Hausdorff distance is related but not
equivalent to local roughness; see \cite{Al3}.  The very-low-temperature
restriction was removed by Ioffe and Schonmann \cite{IS}, who proved
Dobrushin-Kotecky-Schlosman theorem up to the critical temperature. 
For Bernoulli percolation the Wulff construction was justified in
\cite{ACC}, and for the FK model this was done in \cite{Al3}.  For
these models the surface tension is given by the inverse of the
exponential rate of decay of the dual connectivity.

Boundary roughness
has been a topic of considerable interest in the physics literature
(see e.g. \cite{KS}).  The heuristics for the local roughness of
$\Gamma_0$, described in \cite{Al3}, are related to the
boundary-roughness heuristics for two-dimensional growth models such
as first-passage percolation that are believed to be governed by the
``KPZ'' theory
(\cite{KPZ}, \cite{LNP}, \cite{NP}), to polymers in two-dimensional
random environments \cite{Pi}, and, as noted in
\cite{Al3}, to the heuristics of rigorously proved results on
longest increasing subsequences of random permutations \cite{BDJ},
which in turn are related to the fluctuations of eigenvalues of
random matrices (see \cite{Jo}).  In all cases for an object of linear
scale $l$ there is known or believed to be roughness of order
$l^{1/3}$ and a longitudinal correlation length of order
$l^{2/3}$.  In the percolation droplet this correlation length
should appear as the typical separation between adjacent extreme
points of the convex hull of $\Gamma_0$.

In \cite{Al3} the \emph{average local roughness}, denoted
ALR($\Gamma_0)$, for the percolation droplet was defined as the area
between the droplet and its convex hull boundary, divided by the
Euclidean length of the convex hull boundary.  It was proved there that
with high probability, for a droplet conditioned to have area at least
$l^2$, the ALR($\Gamma_0)$ is
$O(l^{1/3} (\log l)^{2/3})$.  The main feature of interest is the
exponent 1/3 matching the KPZ heuristic; the power of $\log l$ may
be considered an artifact of the proof.  Here we consider not
average but \emph{maximum local roughness}, denoted
MLR($\Gamma_0)$ and defined as the maximum distance from any point of
$\Gamma_0$ to the convex hull boundary, and we show that for the
Bernoulli percolation droplet, for some $c_0 > 0$, with high probability
it is at least $c_0 l^{1/3} (\log l)^{-2/3}$.  It was proved in
\cite{Al3} that with high probability 
MLR($\Gamma_0)$ is $O(l^{2/3}(\log l)^{1/3})$, but this is a
presumably a very crude bound, lacking the right power of $l$;
it is more reasonable to compare the lower bound here on
MLR($\Gamma_0)$ to the upper bound for ALR($\Gamma_0)$, as the two
should differ by at most a
multiplicative factor that is a power of
$\log l$, as we explain next.

One way to obtain more-detailed heuristics for the droplet boundary
is to view it as having Gaussian fluctuations about a fixed Wulff
shape of area $l^2$, a point of view justified in part by the
results in \cite{DH} and \cite{Hr}.  This point of view suggests that
if we take a Brownian bridge on [0,1], rescale it by $2\pi l$
horizontally and
$l^{1/2}$ vertically, and wrap it around a circle of radius $l$,
joining
$(0,0)$ and $(2\pi l,0)$, the result should resemble the droplet
boundary. In \cite{Uz} it was proved that for this wrapped Brownian
bridge the maximum local roughness is with high probability bounded
between
$c_1 l^{1/3} (\log l)^{2/3}$ and $c_2 l^{1/3} (\log l)^{2/3}$ for
some $0 < c_1 < c_2 < \infty$.  The exponent 2/3 on $\log l$ here is
related to the L\'evy modulus of continuity for Brownian motion.  The
wrapped-Brownian-bridge heuristic suggests that ALR($\Gamma_0)$
should be of order $l^{1/3}$, without a power of $\log l$, supporting
the idea that ALR($\Gamma_0)$ and MLR($\Gamma_0)$ differ by only a
multiplicative factor that is roughly a power of $\log l$.  The circle 
provides a reasonable heuristic here because Ioffe and Schonmann
\cite{IS} showed that for fixed $p$ the curvature of the boundary of the
unit-area Wulff  shape is bounded away from 0 and $\infty$.

\section{Definitions, Preliminaries, Statement of Main Result}

\par A $bond$, denoted $\langle xy \rangle$, is an unordered pair of
nearest neighbor sites $x,y \in \mathbb{Z}^{2}$. The set of all bonds
between the nearest neighbor sites of $\mathbb{Z}^{2}$, will be
denoted by $\mathbb{B}_{2}$.  Let  $\{ \omega(b), b \in
\mathbb{B}_2\}$ be an i.i.d. family of Bernoulli random variables
with $P(\omega(b) = 1) = p$. Given a realization of $\omega$, a bond
$b \in
\mathbb{B}_{2}$ is said to be $open$ if $\omega(b)=1$ and $closed$
if $\omega(b)=0$. Consider the random graph containing the
vertex set of $\mathbb{Z}^{2}$ and the open bonds only; the
connected components of this graph are called \emph{open clusters}.
For $p$ below the critical probability
$p_{c} = 1/2$ \cite{Ke} all open clusters are
finite with probability one and when
$p > p_{c}$, there exists a unique infinite cluster of open
bonds with probability one. 

\par For $x \in \mathbb{Z}^2$ let $x^*$ denote $x + (1/2,1/2)$.
The lattice with vertex set $\{x^*: \ x\in
\mathbb{Z}^{2}\}$ and all nearest neighbor bonds is called the $dual$
$lattice$.  Each bond $b$ has a unique dual bond, denoted $b^*$, which
is its perpendicular bisector; $b^*$ is defined to be open precisely
when $b$ is closed, so that the dual configuration is Bernoulli
percolation at density $1-p$.
A $(dual)$ $path$ is a sequence $(x_{0}, \langle
x_{0}x_{1} \rangle, x_{1}, \cdots,$
$\langle x_{n-1}, x_{n} \rangle)$ of alternating (dual) sites and
bonds.  A \emph{(dual) circuit} is a path with $x_{n}=x_{0}$ which
has all bonds distinct and does not cross itself (in the obvious
sense).  Note we allow a circuit to touch itself without
crossing, i.e. nondistinct sites are not restricted to $x_n = x_0$.
For a (dual) circuit $\gamma$, the interior
$\Int(\gamma)$ is the union of the bounded components of
the complement of $\gamma$ in $\mathbb{R}^{2}$. An open dual circuit
$\gamma$ is called an \emph{exterior dual circuit} in a configuration
$\omega$ if $\gamma \ \cup \Int(\gamma)$ is maximal among all open dual
circuits in $\omega$. A site $x$ is surrounded by at most one
exterior dual circuit; when this circuit exits, it is denoted by
$\Gamma_{x}$. $|\cdot|$ denotes the Euclidean norm for vectors,
cardinality for finite sets and Lebesgue measure for regions in
$\mathbb{R}^{2}$, depending on the context. For $x,y\in
\mathbb{R}^{2}$, let dist$(\cdot,\cdot)$ and diam$(\cdot)$ denote
Euclidean distance and Euclidean diameter, respectively. Let
$B_{r}(x)$, denote the open Euclidean ball of radius $r$ about
$x$.
 For
$A,B \subset \mathbb{R}^{2}$, define
dist$(A,B)=\inf\{\dist(x,y):x\in A, y \in B\}$ and
dist$(x,A)=$dist$(\{x\},A)$.  We define the
\emph{average local roughness} of a circuit $\gamma$ by
\begin{eqnarray*}
  \mbox{ALR}(\gamma)=\frac{|\mbox{Co}(\gamma) \setminus
  \Int(\gamma)|}{|\partial\mbox{Co}(\gamma)|},
\end{eqnarray*}
where Co$(\cdot)$ denotes the convex hull. The $maximum$ $local$
$roughness$ is 
\begin{eqnarray*}
  \MLR(\gamma)=\sup\{\dist(x,\partial\mbox{Co}(\gamma))\
  : \ x \in \gamma \}
\end{eqnarray*}

Throughout the paper, $K_1,K_2,...$ represent constants
which depend only on $p$.  
Our main result is the following.

\begin{Theorem}\label{E: mainthm}
Let $1/2 < p < 1$.
There exists $K_{1}>0$ such that, under the measure
\newline
$P\big(\ \cdot \ \big| \ | \Int(\Gamma_{0})| \geq
l^{2}\big)$, with probability approaching $1$ as $l\to \infty$ we
have
\begin{eqnarray} \label{E:mainresult}
  \MLR(\Gamma_{0})\geq K_{1}l^{1/3}(\log l)^{-2/3}
\end{eqnarray}
\end{Theorem}

The main ingredients of the proof will be coarse graining
concepts, the renewal structure of long dual connections in the
supercritical regime and exchangeability of the increments 
between regeneration points, all of which will be discussed below.   The
basic idea is that MLR($\Gamma_{0}) < K_{1}l^{1/3} (\log
l)^{-2/3}$ implies that $\Gamma_{0}$ stays in a narrow tube along
its own convex hull, which is a highly unlikely event, due to
the Gaussian fluctuations of connectivities.  More precisely, if $w$
and $w'$ are extreme points of $\Co(\Gamma_0)$ separated by a
distance of order $l^{2/3}(\log l)^{-1/3}$, then MLR($\Gamma_{0}) < 
K_{1}l^{1/3} (\log
l)^{-2/3}$ requres that $\Gamma_{0}$ stay confined within $O(l^{1/3}
(\log l)^{-2/3})$ of the straight line from $w$ to $w'$. 
Gaussian fluctuations, though, would say that the typical deviation
from the straight line is of order $l^{1/3}(\log l)^{-1/6}$, which is the
square root of the length of the line.  Thus the confinement for the
segment between
$w$ and $w'$ is analogous to keeping the maximum magnitude of a
Brownian bridge below $O((\log l)^{-1/2})$, and such confinement along
the entire boundary of $\Gamma_0$ is very unlikely.  The Brownian
bridge analogy is an underlying heuristic but does not enter directly
into our proofs.

\par We use some notation, results and techniques introduced in
\cite{Al3}.  For a family of bond percolation models
including Bernoulli percolation and the FK model, upper bounds have been
established in [Al3] for ALR$(\Gamma_{0})$, MLR($\Gamma_{0})$ and
the deviation between $\partial \Gamma_{0}$ and Wulff shape. We
denote the unit Wulff shape (i.e. the set $W$ of (\ref{E:
classicalW}), normalized to have area 1) by
$\mathbf{K}_{1}$. There
exists constants
$K_{i}$ such that the following hold with probability
approaching to
$1$, as $l
\to \infty$, under the measure $P(\cdot \ | \
| \Int(\Gamma_{0})| \geq l^{2})$:
\begin{eqnarray}
  \mbox{ALR}(\Gamma_{0}) \leq K_{2}l^{1/3} (\log l)^{2/3}, 
  \label{E:ALRbound}\\
  \inf_{x}\ \dist_{H}\big(\partial \mbox{Co} (\Gamma_{0}), x+
  \partial(l
  \mathbf{K}_{1})\big)
  \leq K_{3} l^{2/3} (\log l)^{1/3}, \\
  \MLR(\Gamma_{0}) \leq K_{4} l^{2/3} (\log l)^{1/3}.
\end{eqnarray}
Together, (\ref{E:mainresult}) and (\ref{E:ALRbound}) suggest that local
roughness is of order
$l^{1/3}$, up to a possible logarithmic correction factor, for
sufficiently large
$l$.

\par We will use two standard inequalities for percolation: the
Harris-FKG inequality \cite{Ha} and the BK
inequality \cite{vdBK}. Let $\mathbb{D} \subset \mathbb{B}_2$ and
$\omega,
\widetilde{\omega}
\in \{0,1\}^{\mathbb{D}}$.  We write
$\widetilde{\omega} \geq \omega$ if all
open bonds in $\omega$ are also open in $\widetilde{\omega}$. 
An event $A
\subset \{0,1\}^{\mathbb{D}}$ is \emph{increasing
(decreasing)} if its indicator function
$\delta_{A}$ is nondecreasing (nonincreasing) according to this partial
order.

\emph{Harris-FKG inequality}.  For Bernoulli percolation, if
$A_{1},A_{2},\cdots, A_{n}$ are all increasing, or all decreasing,
events, then
\begin{eqnarray*}
P(A_{1} \cap A_{2}\cap \cdots \cap  A_{n}) \geq
P(A_{1})P(A_{2})\cdots P(A_{n}).
\end{eqnarray*}

For sets $S \subset \mathbb{B}_2$, we will denote by $\omega_{S}$
the restriction of $\omega$ to $S$. 
The event $A$ is said to \emph{occur on} the
set $S$ in the configuration $\omega$ if $\omega'_S =
\omega_S$ implies $\omega' \in A$.  Two events
$A_{1}$ and
$A_{2}$ \emph{occur disjointly} in $\omega$, denoted by
$A_{1}
\circ A_{2}$, if there exist disjoint sets $S_1, S_2$ (depending on
$\omega$) such that
$A_1$ occurs on $S_1$, and $A_2$ occurs on $S_2$, in $\omega$.  The
event that $A_1$ and $A_2$ occur disjointly is denoted $A_1 \circ A_2$.

\emph{BK inequality}.  If $A_{1},\cdots,A_{n}$ are all increasing,
or all decreasing, events then
\begin{eqnarray*}
P(A_{1} \circ A_{2} \circ \cdots \circ A_{n}) \leq 
P(A_{1})P(A_{2})\cdots P(A_{n}).
\end{eqnarray*}

\par Two points $x,y \in (\mathbb{Z}^{2})^*$
are connected, an event written $\{ \ x \longleftrightarrow y\ \}$, if
there exists a path of open dual bonds leading from $x$ to $y$.  The
Harris-FKG inequality implies that $-\log P(0 \leftrightarrow x)$ is a
subadditive function of $x$, and therefore the limit
\begin{eqnarray*}
\tau(x)=\lim_{n\to\infty}-\frac{1}{n}\log P(0^*\leftrightarrow
(nx)^*),
\end{eqnarray*}
exists for $x\in \mathbb{Q}^{2}$, where the limit is taken through
the values of $n$ satisfying $nx\in\mathbb{Z}^{2}$.
This definition extends to $\mathbb{R}^{2}$ by continuity (see
\cite{ACC}). $\tau$ is a strictly convex norm on $\mathbb{R}^{2}$; the
strict convexity is shown in \cite{CI}.
The $\tau$-norm for unit vectors serves as the surface tension for our
context. Let $\mathbb{S}$ denote the unit
circle in $\mathbb{R}^{2}$. It is known (\cite{Al2},\cite{Me}) that for
$1/2 < p < 1$,
\begin{equation}
0< \min_{x\in \mathbb{S}} \tau(x) \leq \max_{x\in \mathbb{S}} \tau(x) <
\infty,
\end{equation}
\begin{eqnarray} \label{E:coniq}
\beta_{1} |x|^{-\beta_{2}}\exp(-\tau(x))\leq P(0^*\leftrightarrow
x^*) \leq \exp(-\tau(x))
\end{eqnarray}
for some constants $\beta_{1},\beta_{2}>0$  and
\begin{eqnarray} \label{E:taubounds}
\frac{\tau(e)}{\sqrt{2}} \leq \frac{\tau(x)}{|x|}\leq \sqrt{2}
\ \tau(e),
\end{eqnarray}
 where $e$ is a coordinate vector.  

For $x,y\in \mathbb{R}^{2}$, let dist$_{\tau}(\cdot,\cdot)$ and
diam$_{\tau}(\cdot)$ denote the $\tau$-distance and the $\tau$-
diameter, respectively.
 Some of the properties of connectivities and
geometry of Wulff shapes  will be given next.
Denote the unit $\tau$-unit ball by $\mathbf{U}_{1}$:
\begin{equation*}
  \mathbf{U}_{1} {=} \bigl\{x \in \mathbb{R}^{2} \ : \
  \tau(x) \leq 1 \bigl\}
\end{equation*}
and the Wulff shape by $\mathbf{W}_{1}$:
\begin{equation*}
  \mathbf{W}_{1} {=} \bigl\{ t \in
  \mathbb{R}^{2}\ : \ (t,z)_{2} \leq \tau(z) \mbox{ for all } z
  \in \mathbb{S} \bigr\},
\end{equation*}
so that $0 \in \Int(\mathbf{W}_{1})$ and $\mathbf{K}_1 =
\mathbf{W}_1 / |\mathbf{W}_1|$.  We also refer to multiples of
$\mathbf{W}_{1}$ as Wulff shapes.  For the functional
\begin{eqnarray*}
  \mathcal{W}(\gamma)=\int_{\gamma}\tau(v_{x})\ dx,
\end{eqnarray*}
$\mathbf{K}_1$ minimizes 
$\mathcal{W}(\partial V)$ over all regions $V$ with piecewise
$C^1$ boundary, subject to the constraint $|V|=1$; here $v_{x}$ is the
unit forward tangent vector at $x$ and
$dx$ is arc length.  (A class
larger than the regions with piecewise $C^1$ boundary can be used here,
but is not relevant for our purposes; for specifics see \cite{Ta1},
\cite{Ta2}.)  We
define the Wulff constant $\mathcal{W}_1 =
\mathcal{W}(\partial \mathbf{K}_{1})$.
For every $t \in \partial
\mathbf{W}_{1}$ and $x\in \partial \mathbf{U}_{1}$, we have
\begin{equation*}
  1= \max_{y\in\mathbf{U}_{1}}(t,y)_{2}= \max_{s\in \partial
  \mathbf{W}_{1}} (s,x)_{2}.
\end{equation*}

\begin{Definition} \label{D:polar}
Given $x\in \mathbb{R}^{2}\backslash \{0\}$, a point $t\in \partial
\mathbf{W}_{1}$ is \emph{polar} to $x$ if
\begin{equation*}
  (t,x)_{2}= \tau(x) =\max_{s\in \partial \mathbf{W}_{1}}(s,x)_{2}
\end{equation*}
\end{Definition}

\section{Renewal Structure of Connectivities}
For the remainder of the paper we assume we have fixed $1/2 < p < 1$.

 This section will  follow  Section 4 of \cite{CI}.
For $x,y \in (\mathbb{Z}^{2})^*$ and  $t \in \partial \mathbf{W}_{1}$,
we define the line
\begin{equation*}
  \mathcal{H}_{x}^{t} = \{ z \in \mathbb{R}^{2} \ | \ (
  t,z )_{2} =  (t,x)_{2} \}
\end{equation*}
and the slab
\begin{equation*}
  \mathcal{S}_{x,y}^{t} = \{ z \in \mathbb{R}^{2} 
  \ \arrowvert (t,x)_{2} \
  \leq (t,z)_{2} \leq (t,y)_{2}\}.
\end{equation*}
When $x$ and $y$ are connected in the restriction of
the percolation configuration to the slab $\mathcal{S}_{x,y}^{t}$
(excluding the bonds that are only partially in
$\mathcal{S}_{x,y}^{t}$), $\mathbf{C}_{x,y}^{t}$ denotes the
set of sites in the corresponding common cluster inside
$\mathcal{S}_{x,y}^{t}$. Let $e = e(t)$ be a unit vector in the direction
of one of the axes such that the scalar product of $e$ with $t$ is
maximal.

\begin{Definition} For $x,y \in (\mathbb{Z}^{2})^*$ 
satisfying $(t,x)_{2}
< (t,y)_{2}$, let $\bigl\{ \ x
\overset{\widetilde{h_{t}}}{\longleftrightarrow} y \ \bigl\}$ denote
the event that $x$ and $y$ are
$\widetilde{h_t}$-\emph{connected}, meaning $x$ and $y$
are connected by an open dual path in $\mathcal{S}_{x,y}^{t}$.
Let $\bigl\{ \ x \overset{h_{t}}{\longleftrightarrow} y \ \bigl\}$
denote the event that $x$ and $y$ are $h_t$-\emph{connected}, 
meaning $x$ and $y$ are connected inside
$\mathcal{S}_{x,y}^{t}$ and
\begin{equation*} 
 \mathbf{C}_{x,y}^{t} \cap \mathcal{S}_{x,x+e}^{t} = \{x,x+e\} \
 \mbox{ and } \
 \mathbf{C}_{x,y}^{t}  \cap \mathcal{S}_{y-e,y}^{t} = \{y-e,y\}.
\end{equation*}
Let $\bigl\{ \ x \overset{f_{t}}{\longleftrightarrow} y \ \bigl\}$
denote the event that $x$ and $y$ are $f_t$-\emph{connected},
meaning $x \overset{h_{t}}{\longleftrightarrow} y$ and for no $z \in
\Int(\mathcal{S}^{t}_{x,y})$ do both $x
\overset{h_{t}}{\longleftrightarrow} z$ and
$z \overset{h_{t}}{\longleftrightarrow} y$.
\end{Definition}

\begin{Definition} Given a configuration and given $x,y$ with $x
\leftrightarrow y$, we say that
$z \in (\mathbb{Z}^{2})^*$ is a \emph{regeneration point} if 
$(t,x)_{2} < (t,z)_{2} < (t,y)_{2}$ and $\mathbf{C}_{x,y}^{t}
\cap \mathcal{S}_{z-e,z+e}^{t}=\{z-e,z,z+e\}$.
\end{Definition}

\par Let $\mathcal{R}_{x,y}^{t}$ denote the random set of regeneration
points of $\mathbf{C}_{x,y}^{t}$. Next, a probabilistic bound on
the size of $\mathcal{R}_{x,y}^{t}$ will be given. For our purposes, we
need a different formulation of Lemma 4.1 of \cite{CI}: we use $\bigl\{
\ x \overset{\widetilde{h_{t}}}{\longleftrightarrow} y \ \bigl\}$
instead of $\bigl\{ \ x \overset{h_{t}}{\longleftrightarrow} y \
\bigl\}$ to state the lemma, but the proof is same with
minor changes.

\begin{Lemma} \label{L:fewrenewals}
For every $\epsilon \in (0,\frac{1}{2})$, there exists
$\lambda>0$, $\delta
>0$ and $\nu >0$ such that for all $t_{0}\in \partial \mathbf{W}_{1}$, 
$t \in B_{\lambda}(t_{0})$ and all
$x$ satisfying $(t,x)_{2}\geq
(1-\epsilon)\tau(x)$ we have
\begin{equation} \label{L: noofrenewals}
P\bigl(|\mathcal{R}_{0,x}^{t_{0}}|< \delta|x| \ ; \ \ 0
\overset{\widetilde{h_{t_{0}}}}{\longleftrightarrow} x \bigl) \leq
\exp\{-(t,x)_{2}-\nu|x|\}.
\end{equation}
\end{Lemma}

\section{Coarse Graining and Related Preliminaries}
We will use the coarse graining setup and results of [Al3].  For $s>0$,
and any contour with a
$\tau$-diameter of at least $2s$, the coarse graining algorithm
selects a subset $\{w_{0},w_1,\cdots,w_{m+1}\}$ of the extreme
points of Co$(\gamma)$,
with $w_{m+1} = w_0$, called the $s$-\emph{hull skeleton} of $\gamma$
and denoted HSkel$_s(\gamma)$. The points $w_i$ of HSkel$_s(\gamma)$
appear in order as one traces $\gamma$ in the direction of positive
orientation. We denote the polygonal path
$w_{0}\to w_1 \to \cdots
\to w_{m+1}$ by HPath$_{s}(\gamma)$.  The specifics of the algorithm
for choosing the $s$-hull skeleton are not important to us here; we
refer the reader to \cite{Al3}.  What we need are the following
properties, also from \cite{Al3}.

\begin{Lemma}  There exist constants
$K_{5},K_{6},K_{7},K_{8}>0$ such that for every $s>0$ and every circuit
$\gamma$ having
$\tau$-diameter at least $2s$, the $s$-hull skeleton
$\HSkel_s(\gamma) = \{w_{0},w_{1},\cdots,w_{m+1} \}$ satisfies
\begin{equation}
  m+1 < \frac{K_{5}\diam(\gamma)}{s} \label{E:cg1},
\end{equation}
\begin{equation}
  |\Int(\gamma) \setminus \Int(\HPath_{s}(\gamma))| \leq
  K_{6}s^{2}, \label{E:cg2}
\end{equation}
\begin{equation}
  \sup_{x \in \mbox{Co}(\gamma)} \dist(x,
  \Int(\HPath_{s}(\gamma)) \leq \frac{K_{7} s^{2}}
  {\diam(\gamma)}, \label{E:skeltohull}
\end{equation}
\begin{equation}
  \mathcal{W}(\partial \Co(\gamma)) \leq
  \mathcal{W}(\HPath_{s}(\gamma))+ \frac{K_{8} s^{2}}
  {\diam(\gamma)}.
\label{E:cg4}
\end{equation}
\end{Lemma}

For $0 < \theta  < 1$ a small constant to be specified later, 
our choice of $s$ is
\begin{eqnarray*} 
  s=\bigg(\frac{\theta \sqrt{\pi}}{2K_{7}}\bigg)^{1/2} l^{2/3}(\log
  l)^{-1/3}.
\end{eqnarray*}
Suppose
HSkel$_{s}(\Gamma_{0}) = \{ w_{0},w_{1}, \dots ,w_{m+1} \}$ with
$w_{m+1}=w_{0}$. We define
\begin{eqnarray*}
  \mathcal{L}= \bigg\{ i : |w_{i+1}-w_{i}|\geq \frac{s\sqrt{\pi}}{16
  K_{5}}\bigg\}
\end{eqnarray*}
For $i\in \mathcal{L}$, we call the side between $w_{i}$ and
$w_{i+1}$ \emph{long}. The next lemma gives a lower bound on the
sum of the lengths of long sides when diam$(\Gamma_{0})$ is
not abnormally large.  From \cite{Al3}, for some
$K_{9}, K_{10}, K_{11} > 0$, for $T > 0$,
\begin{eqnarray*} 
  P(\diam_{\tau}(\Gamma_{0}) \geq T) \leq K_{9} T^{4}e^{-T}
\end{eqnarray*}
and
\begin{equation} \label{E:droplowerbd}
  P\bigl( |\Int(\Gamma_0)| \geq l^2 \bigr) \geq
  K_{10} \exp \big(-\mathcal{W}_1 l - K_{11} l^{1/3}(\log l)^{2/3} \big),
\end{equation}
so that for large $l$,
\begin{equation*} 
  P \bigl( \diam_{\tau}(\Gamma_{0}) \geq 2 \mathcal{W}_1 l \ \big| \
  |\Int(\Gamma_0)| \geq l^2 \bigr) \leq e^{-\mathcal{W}_1 l/2}.
\end{equation*}
Also using $\eqref{E:taubounds}$, we have
\begin{eqnarray*}
  \diam(\Gamma_{0}) \leq \frac{\sqrt{2}}{\tau(e)} \
  \diam_{\tau}(\Gamma_{0}) \leq \frac{4 \sqrt{2}}{\mathcal{W}_1}
  \diam_{\tau}(\Gamma_{0})
\end{eqnarray*}
 where in the second inequality we use
$\mathcal{W}_1 \leq 4 \tau(e)$, which follows from the
fact that the unit square encloses the unit area. Therefore
\begin{equation} \label{E:Eucldiam}
  P \bigl( \diam(\Gamma_{0}) \geq 8 \sqrt{2} l \ \big| \
  |\Int(\Gamma_0)| \geq l^2 \bigr) \leq e^{-\mathcal{W}_1 l/2}.
\end{equation}
so to prove
Theorem \ref{E: mainthm} we need only consider configurations with
diam$(\Gamma_{0})
< 8\sqrt{2}\ l$.  We say that $\{w_0,..,w_{m+1}\}$ is
$l$-\emph{regular} if there exists a configuration in which
$|\Int(\Gamma_0)| \geq l^2$, diam$(\Gamma_0) < 8 \sqrt{2}\ l$
and HSkel$_s(\Gamma_0) = \{w_0,..,w_{m+1}\}$.

\begin{Lemma} \label{L:Lsum}
If $\{
w_{0},w_{1}, \dots ,w_{m+1} \}$ is $l$-regular and $l$ is sufficiently
large, then
\begin{eqnarray} \label{E:sumlong}
\sum_{i \in \mathcal{L}}|w_{i+1}-w_{i}| \geq
\sqrt{\frac{\pi}{2}} \ l
\end{eqnarray}
\end{Lemma}

\begin{proof}
$\eqref{E:cg2}$ implies that for some $K_{12}$, and $\Gamma_0$ as in the
definition of $l$-regular,
\begin{eqnarray*}
|\mbox{Int(HPath}_{s}(\Gamma_0))| \geq l^{2}- K_{12}l^{4/3}(\log
l)^{-2/3}
\geq \frac{l^{2}}{2}\ ,
\end{eqnarray*}
 where the last inequality is satisfied for
sufficiently large $l$. By the standard isoperimetric inequality,
it follows that
\begin{eqnarray*}
\sum_{i \in \mathcal{L}}|w_{i+1}-w_{i}|+ \sum_{i \in
\mathcal{L}^{c}}|w_{i+1}-w_{i}| \geq l \sqrt{2\pi}.
\end{eqnarray*}
Using $\eqref{E:cg1}$, the total number of sides can be bounded above:
\begin{eqnarray*}
m+1 \leq \frac{K_{5} \diam(\Gamma_{0})}{s} \leq  \frac{8
\sqrt{2} \ K_{5} \ l}{s}.
\end{eqnarray*}
Therefore
\begin{eqnarray*}
 \sum_{i
\in \mathcal{L}^{c}}|w_{i+1}-w_{i}| \leq (m+1)
\frac{s\sqrt{\pi}}{16K_{5}}  \leq  \sqrt{\frac{\pi}{2}} l,
\end{eqnarray*}
and the lemma follows.
\end{proof}

We next need to specify the vector $t_i$ which will be used to
define slabs and regeneration points for the connection from $w_i$ to
$w_{i+1}$.  The natural choice is to take
$t_{i}$ polar to $w_{i+1}-w_{i}$, but
in order to avoid some technicalities in upcoming proofs we
will choose $t_{i}$ to be close to the polar value, but having
rational slope.  Let $V \subset \mathbb{R}^2$ denote the wedge
consisting of those vectors $x$ such that the 
angle from the positive horizontal axis to $x$ is
in $[0,\pi/4]$.  Due to lattice symmetries we may assume that
$w_{i+1}-w_{i} \in V$.  Let $\widetilde{t_{i}} \in
\partial \mathbf{K}_{1} \cap V$ be such that $\widetilde{t_{i}}$ is
polar to $w_{i+1}-w_{i}$.
Then the angular difference
between $\widetilde{t_{i}}$ and $w_{i+1}-w_{i}$  is at most
$\pi/4$.  The existence of a polar point with such
properties is guaranteed by symmetries of $\mathbf{K}_{1}$. Let us fix
$\epsilon\in(0,1/2)$, and let
$\lambda=\lambda(\epsilon)\mbox{ as in } \eqref{L: noofrenewals}$.
We choose $t_{i} \in V \cap B_{\lambda}(\widetilde{t_{i}}) \cap
\partial
\mathbf{K}_{1}$ so that the slope of $t_i$ is $r/q$, with
$q=[1/\lambda]+1$ and $r \in \mathbf{Z}$.  Choosing $t_{i}$ this way
will allow us to use$~\eqref{L: noofrenewals}$, with the parameters
$t_{0}$ and
$t$ chosen as $t_{i}$ and $\widetilde{t_{i}}$,
respectively.  Note that $e(t_i) = (1,0)$, which we denote by $e_i$.

\input{epsf}
\begin{figure} 
\epsfxsize=6.0in
\begin{center}
\leavevmode
\epsffile{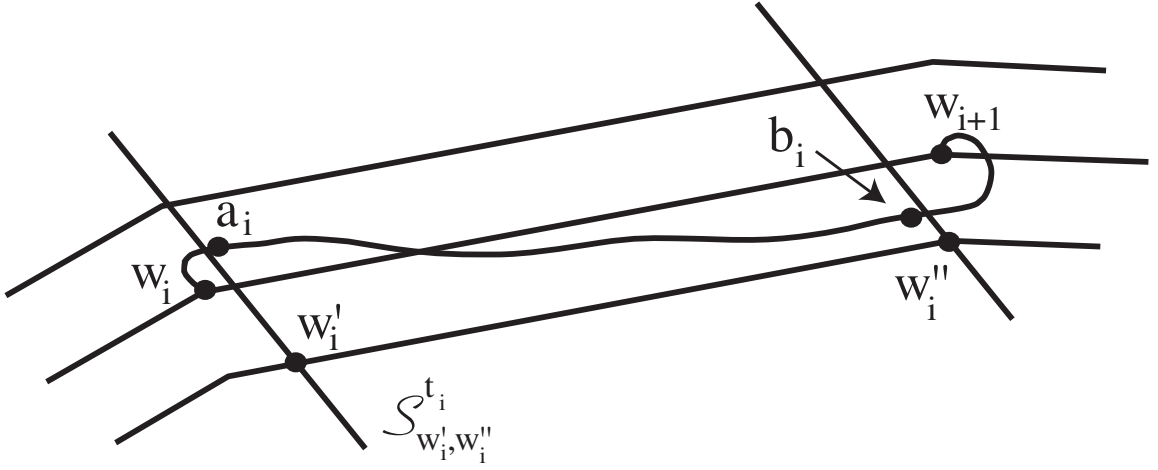}
\end{center}
\caption{A section of $A_d$, and a connection from $w_i$ to $w_{i+1}$
which includes a cylinder connection from $a_i$ to $b_i$.}
\label{F:Adsection}
\end{figure}

\par By $\eqref{E:skeltohull}$ for our chosen $s$, the deviation
between 
$\Co(\Gamma_{0})$ and Int(HPath$_{s}(\Gamma_{0}))$ inside
it does not exceed
$\theta l^{1/3}(\log l)^{-2/3}$. Let $l_{i}$ be the line through
$w_{i}$ and $w_{i+1}$.  We set $d= 2 \theta l^{1/3}(\log
l)^{-2/3}$, and define $A_{d}$, the annular tube of diameter $2d$
around HSkel$_{s}(\Gamma_{0})$, as follows.
Denote the line parallel to $l_{i}$ which is $d$ units
outside of HSkel$_{s}(\Gamma_{0})$ by $l_{i}^{+}$ and
the line parallel to $l_{i}$ which is $d$ units in the opposite
direction by $l_{i}^{-}$. Let
$H_{l_i}$ be the half space bounded by 
$l_{i}$ that contains 
HSkel$_{s}(\Gamma_{0})$, let $H_{l_i^{\pm}}$ be the
halfspaces bounded by $l_i^{\pm}$ such that $H_{l_i^-} \subset
H_{l_i} \subset H_{l_i^+}$ and let
\begin{eqnarray*}
  A_{d}= A_d(w_0,..,w_{m+1}) = 
  \biggl(\bigcap_{i=1}^{m} \ H_{l_{i}^{+}} \biggl) \ \setminus \
  \biggl( \bigcap_{i=1}^{m}
  \ H_{l_{i}^{-}} \biggl)
\end{eqnarray*}
(see Figure \ref{F:Adsection}.)
Let $T_{d}^{i}$ denote the (infinite) tube with diameter $2d$,
bounded by $l_{i}^{+}$ and $l_{i}^{-}$. 
Let $w'_{i}$ and $w''_i$ be the points on
$l_{i}^{-}$ such that $\mathcal{S}^{t_{i}}_{w'_{i},w''_i}$ is
the largest slab satisfying
\begin{eqnarray*}
  \mathcal{S}^{t_{i}}_{w'_{i},w''_i} \cap T_{d}^{i} \cap A_d =
  \mathcal{S}^{t_{i}}_{w'_{i},w''_i} \cap \ T_{d}^{i}.
\end{eqnarray*} 
Let $B_i$ be the event that
there exist $a_{i} \in
\mathcal{S}_{w'_{i},w'_{i}+e_i}^{t_{i}} \cap T_{d}^{i}$ and $ \ b_{i} \in
\mathcal{S}_{w''_i-e_i,w''_i}^{t_{i}} \cap T_{d}^{i}$ such that
the event
\begin{eqnarray*} 
  \ \ \{w_{i}
  \longleftrightarrow a_{i}\}\ \circ \
  \{a_{i} \overset{\widetilde{h_{t_i}}}{\longleftrightarrow} b_{i}
  \mbox{ in } T^i_{d}\}\ \circ \{\ b_i
  \longleftrightarrow w_{i+1}\}
\end{eqnarray*}
occurs.  For configurations in $\big\{ w_i \longleftrightarrow
w_{i+1}$ in $A_d
\big\} \backslash B_i$, every open path from $w_i$ to $w_{i+1}$ must
go ``the long way around $A_d$''; presuming $l$ is
large and $\{w_0,..,w_{m+1}\}$ is
$l$-regular, for some $K_{13}$ this implies that $w_i \leftrightarrow
z$ for some $z \in
\mathcal{S}_{w_{i},w_{i}+e_i}^{t}$ with dist$(z,w_i) \geq K_{13}l$.  By
(\cite{Al3}, Lemma 7.1) we then have for some $K_{14},K_{15}$,
\begin{equation} \label{E:noBibound}
  P\bigl( B_i^c \ \big| \ w_i \longleftrightarrow
  w_{i+1} \bigr) \leq K_{14}e^{-K_{15}l}.
\end{equation}

\begin{Lemma} There exists constants $K_{14},K_{15}>0$
such that for $\{w_0,..,w_{m+1}\}$ $l$-regular and $\epsilon, t_i$ as in
the preceeding, we have
\begin{align} \label{E: Adtotubelemma}
  P( w_{i} &\leftrightarrow w_{i+1} \mbox{ in }A_{d}\ \big| w_{i}
    \leftrightarrow w_{i+1}) \\
  &\leq K_{14}\exp{(-K_{15}l)+\sum_{a_{i},b_{i}}}
    \ P( a_{i} \leftrightarrow b_{i}
    \mbox{ in } T_{d}^{i} \  \bigl| \    a_{i}
    \overset{\widetilde{h_{t_{i}}}}{\leftrightarrow} b_{i}) \notag
\end{align}
where the sum is over all $a_{i} \in
\mathcal{S}^{t_{i}}_{w'_{i},w'_{i}+e_i} \cap T_{d}^{i} \cap
(\mathbb{Z}^2)^*$  and
$b_{i}\in \mathcal{S}^{t_{i}}_{w''_i, w''_i-e_i} \cap T_{d}^{i}
\cap (\mathbb{Z}^2)^*$.
\end{Lemma}

\begin{proof}
By (\ref{E:noBibound}) we can bound $P( w_{i} \leftrightarrow
w_{i+1} \mbox{ in }A_{d})$ by
\begin{equation}
  K_{14}e^{-K_{15}l} P(w_{i} \leftrightarrow w_{i+1})
  +\sum_{a_{i},b_{i}} P\big(\{w_{i}
  \leftrightarrow a_{i}\} \circ  \{
  a_{i} \overset{\widetilde{h_{t_i}}}{\longleftrightarrow} b_{i}
  \mbox{ in } T^i_{d}\} \circ  \{b_{i}
  \leftrightarrow w_{i+1}\} \big), \notag
\end{equation}
where the sum is over all $a_{i} \in
\mathcal{S}^{t_{i}}_{w'_{i},w'_{i}+e_i} \cap T_{d}^{i} \cap 
(\mathbb{Z}^2)^*$  and
$b_{i}\in \mathcal{S}^{t_{i}}_{w''_i,w''_i-e_i} \cap T_{d}^{i} \cap
(\mathbb{Z}^2)^*$.  We now apply the BK and FKG inequalities:
\begin{eqnarray}
&\sum_{a_{i},b_{i}}&P\bigg(\{w_{i} \leftrightarrow a_{i}\}\ \circ \
  \{a_{i} \overset{\widetilde{h_{t_i}}}{\leftrightarrow} b_{i} \mbox{
  in } T^i_{d}\}\ \circ \ \{ b_i
  \leftrightarrow w_{i+1}\} \bigg) \notag \\
&\leq &\sum_{a_{i},b_{i}}  P(w_{i} \leftrightarrow a_{i}) \ P( a_{i}
  \overset{\widetilde{h_{t_i}}}{\leftrightarrow} b_{i} \mbox{ in }
  T_{d}^{i}) \ P(b_i \leftrightarrow w_{i+1}), \notag \\
&=&\sum_{a_{i},b_{i}}P(w_{i} \leftrightarrow a_{i}) \
  P(a_{i}\overset{\widetilde{h_{t_i}}}{\leftrightarrow} b_{i}) \ P(
  a_{i} \overset{\widetilde{h_{t_i}}}{\leftrightarrow} b_{i} \mbox{
  in } T_{d}^{i}\ | \
  a_{i}\overset{\widetilde{h_{t_i}}}{\leftrightarrow} b_{i}) \ P(b_i
  \leftrightarrow w_{i+1}) \notag \\
&\leq& \sum_{a_{i},b_{i}} P(w_{i} \leftrightarrow w_{i+1}) \ P(
  a_{i} \overset{\widetilde{h_{t_i}}}{\leftrightarrow} b_{i} \mbox{
  in } T_{d}^{i}\ | \
  a_{i}\overset{\widetilde{h_{t_i}}}{\leftrightarrow} b_{i}), \notag
\end{eqnarray}
and$~\eqref{E:
Adtotubelemma}$ follows.
\end{proof}

In order to bound the probability of the event $\bigl\{a_{i}
\overset{\widetilde{h_{t_i}}}{\longleftrightarrow} b_{i}
\mbox{ in } T_{d}^{i}\bigl\} $ using the
renewal structure of cylinder connectivities, we need control of
the size of
$|b_{i}-a_{i}|$ to apply$~\eqref{L: noofrenewals}$.  The parallelogram 
$\mathcal{S}_{w'_{i}+e,w''_i-e}^{t_{i}}\cap T_{d}^{i}$ has 2 short
sides (the sides not parallel to $w_{i+1} - w_i$), one near $w_i$ and
the other near $w_{i+1}$ (see Figure \ref{F:Adsection}).  It follows
easily from the fact that
$w_{i+1}-w_i, t_i$ are in the wedge $V$ that for every $a$ in the short
side near $w_i$ we have $|w_i - a| \leq 2d\sqrt{2}$, and analogously for
$w_{i+1}$.  Therefore
\[
  |w_i - a_i| \leq 2d\sqrt{2} + 1, \qquad |w_{i+1} - b_i| \leq
  2d\sqrt{2} + 1,
\]
and hence 
\begin{equation} \label{E:wbarelation}
\bigl|(w_{i+1}-w_{i}) - (b_{i}-a_{i})\bigl| \leq   4d\sqrt{2}+2.
\end{equation}
Since
\begin{equation} \label{E:polarprop}
  \tau (w_{i+1}-w_{i})=(\widetilde{t_{i}},w_{i+1}-w_{i})_{2}, 
\end{equation}
provided
$l$ is large we have
\begin{equation} \label{E:baapprox}
  (\widetilde{t_{i}},b_i - a_i)_{2} \geq (1 - \epsilon)\tau(b_i - a_i)
\end{equation}
for our chosen $\epsilon$.

\par
\begin{Lemma}
Given $\epsilon, t_i, a_i, b_i$ as in the preceeding
and $\delta$ as in $\eqref{L:
noofrenewals}$, there exists
$\nu'>0$ such that provided $l$ is sufficiently large,
\begin{eqnarray} \label{L: noofrenewals2}
P\big(|\mathcal{R}_{a_{i},b_{i}}^{t_{i}}|<\delta |b_{i}-a_{i}| \  \big| \
a_{i}
\overset{\widetilde{h_{t_{i}}}}{\leftrightarrow} b_{i} \big) \leq \exp
(-\nu'|b_{i}-a_{i}|).
\end{eqnarray}
\end{Lemma}

\begin{proof}
From (\cite{Al3} equation (7.6)), for some
$K_{16},K_{17}>0$, we have
\begin{eqnarray} \label{E:aibilower}
  {P\big(a_{i} \overset{\widetilde{h_{t_{i}}}}{\leftrightarrow}
  b_{i} \big) } \geq K_{16} |b_{i}-a_{i}|^{-K_{17}}
  \exp\big(-\tau(b_{i}-a_{i})\big).
\end{eqnarray}
By (\ref{E:baapprox}), Lemma \ref{L:fewrenewals} applies; with
(\ref{E:aibilower}) this shows that for some
$\nu > 0$,
\begin{align}
  &P\big(|\mathcal{R}_{a_{i},b_{i}}^{t_{i}}|<\delta |b_{i}-a_{i}|
  \  \big| \ a_{i} \overset{\widetilde{h_{t_{i}}}}{\leftrightarrow}
  b_{i} \big) \label{E:interm} \\ \nonumber
  &\leq
  \frac{1}{K_{16}}|b_{i}-a_{i}|^{K_{17}}\exp\big(-(\widetilde{t_{i}},
  b_{i}-a_{i})_{2}+\tau(b_{i}-a_{i})-\nu
  |b_{i}-a_{i}|\big).
\end{align}
By (\ref{E:wbarelation}) and (\ref{E:polarprop}), we have
\begin{align}
 &-(\widetilde{t_{i}},b_{i}-a_{i})_{2}+\tau(b_{i}-a_{i})-\nu 
  |b_{i}-a_{i}| \notag \\ 
 &\quad \leq 2 \tau(w_{i+1}-w_{i} -b_{i}+a_{i})-\nu|b_{i}-a_{i}| \notag
   \\
 &\quad \leq  K_{18}(4d\sqrt{2}+2)- \nu |b_{i}-a_{i}| \notag
\end{align}
for some $K_{18}>0$. Since $d$ is small compared to
$|b_{i}-a_{i}|$, using this bound in$~\eqref{E:interm}$, for some
constant $\nu'<
\nu$ we have$~\eqref{L: noofrenewals2}$.
\end{proof}

\par
  Next, we will define orthogonal increments between adjacent
regeneration points. There is no canonical choice of direction relative
to which increments are defined; we will use the direction 
orthogonal to the line joining $w_{i}$ and $w_{i+1}$.
\begin{Definition} For any $x \in \mathcal{S}_{w_{i},w_{i+1}}^{t_{i}}$,
define
$f:\mathcal{S}_{w'_{i},w''_i}^{t_{i}} \to \mathbb{R}$ as follows:
\begin{equation*}
  f(x)=\begin{cases}\hfill \dist(x,l_{i}), & \mbox{ if $x$ is above
  the line } l_{i},\mbox{ joining } w_{i} \mbox{ and } w_{i+1}, \\
  -\dist(x,l_{i}), & \mbox{ if $x$ is on or below the line }
  l_{i}.
  \end{cases}
\end{equation*}
\end{Definition}
For the following definitions assume $a_{i} 
\overset{\widetilde{h_{t_{i}}}}{\leftrightarrow}
b_{i}$. The regeneration points between $a_i$ and
$b_i$ have a natural ordering according to their distance from
$\mathcal{H}_{a_i}^{t_i}$.
\begin{Definition}
For  $r' \in \mathcal{S}^{t_{i}}_{a_{i},b_{i}}$ define
$\Delta:\mathcal{S}^{t_{i}}_{a_{i},b_{i}} \to \mathbb{R}$ as follows:
\begin{eqnarray*}
  \Delta(r') =\begin{cases} f(r'), & \mbox{ if } r' 
  \mbox{ is the first regeneration point}, \\
 f(r')-f(\tilde{r}), & \mbox{ if } \tilde{r},r' \mbox{ are
  successive regeneration points},\\
  \quad 0 & \mbox{ if } r' \mbox{ is not a regeneration
  point.} \end{cases}
\end{eqnarray*}
\end{Definition}
\begin{Definition} For $\mathcal{H}_{z}^{t_{i}} \subset  
\mathcal{S}^{t_{i}}_{a_{i},b_{i}}$ define
\begin{eqnarray*}
\widetilde{\Delta}(\mathcal{H}_{z}^{t_{i}}) = \begin{cases} \Delta(r') &
\mbox{ if there is a regeneration point } r' \in
\mathcal{H}_{z}^{t_{i}}, \\ \quad 0 & \mbox{ otherwise.}
\end{cases}
\end{eqnarray*}
\end{Definition}
We will refer to the values $\Delta(r)$ as \emph{increments}. We need to
show that, given $a_{i} \overset{\widetilde{h_{t_{i}}}}{\leftrightarrow}
b_{i}$, there are unlikely to be too many small increments.  This
will be proved by showing that a positive proportion of increments have
magnitude greater than equal to $1/2$, with high probability. This result
will be used to bound the variance of sums of increments from below.

\par
For $\delta$ as in Lemma \ref{L:fewrenewals}, and $a_i,b_i$ fixed, let
$N=\lfloor
\delta|b_{i}-a_{i}| \rfloor$, and $R=\lfloor N/8 \rfloor$.  Let
$\mathcal{U}$ be the collection of all $(z_{1},\cdots,z_{R})$ such that
for $ j=1, \cdots,R$, we have
\begin{enumerate}
\item[(i)] $z_{j}\in
\mathcal{S}^{t_{i}}_{a_{i},b_{i}}$;
$z_{j}$ is on the line through
$w_{i}$, parallel to $t_{i}$,
\item[(ii)] $(t_{i},z_{1})_{2}< (t_{i},z_{2})_{2} < \cdots <
(t_{i},z_{R})_{2}$,
\item[(iii)] $(\Int\mathcal{S}^{t_{i}}_{z_{j}-4e_i,z_{j}+4e_i})$ and
$(\Int\mathcal{S}^{t_{i}}_{z_{k}-4e_i,z_{k}+4e_i})$ 
are disjoint for $j \neq k$.  
\end{enumerate}
By property (i), there is a bijection pairing
$\{z_{1},z_{2},\cdots,z_{R}\}\in \mathcal{U}$ and the set of lines
$\mathcal{H}_{z_{j}}^{t_{i}}$ passing through the points
$\{z_{1},z_{2},\cdots,z_{R}\}$. Suppose
$\mathcal{R}^{t_{i}}_{a_{i},b_{i}}=\{r_{1},r_{2},\cdots,r_{I}\}$, with
$I \geq N$. Next, we define $\mathcal{Q}^{t_{i}}_{a_{i},b_{i}} = \{
\sigma_1,..,\sigma_R \} \subset \mathcal{R}^{t_{i}}_{a_{i},b_{i}}$
according to the following algorithm:  
\begin{enumerate}
\item[(1)] $\sigma_{1}= r_{k_{1}}$, where $k_{1}$ is the
smallest integer satisfying
$(t_{i},a_{i}+4e_i)_{2}\leq(t_{i},r_{k_{1}})_{2}$, 
\item[(2)]
$\sigma_{j}=r_{k_{j}}$, where $k_{j}$ is the smallest integer
satisfying $(t_{i},\sigma_{j-1}+8e_i)_{2}\leq(t_{i},r_{k_{j}})_{2}$,
for $j=2,3,\cdots,R.$ 
\end{enumerate}
For $j\geq 2$, this algorithm can
skip at most 7 regeneration points after $\sigma_{j-1}$ before it
selects
$\sigma_{j}$; under the assumption that there are at least 
$N$ regeneration points, it will successfully choose exactly $R$
regeneration points.  ($\mathcal{Q}^{t_{i}}_{a_{i},b_{i}}$ is undefined
when there are fewer than $N$ regeneration points, so
$|\mathcal{Q}^{t_{i}}_{a_{i},b_{i}}| = R$ whenever
$\mathcal{Q}^{t_{i}}_{a_{i},b_{i}}$ is defined.)  Notice that, for some
$(z_{1},z_{2}, \cdots,z_{R})\in
\mathcal{U}$, the regeneration point $\sigma_{j}$ occurs on
$\mathcal{H}_{z_{j}}^{t_{i}}$, for $j=1,2,\cdots, R$.  Also, since the
slope of $t_i$ is rational, the line $\mathcal{H}_{\sigma_j}^{t_i}$
contains other lattice points, which are also possible locations for the
$j$th regeneration point, when only
$\mathcal{H}_{\sigma_j}^{t_i}$ is specified.   

\begin{Lemma}
Given $\epsilon, t_i, a_i, b_i$ as in the preceeding,
for $\delta>0$ from$~\eqref{L: noofrenewals}$, there exist
$\gamma, \varphi >0$ such that
\begin{equation}
  P\biggl( \sum_{k=2}^{N} \delta_{\{|\Delta (r_{k})| \geq
    \frac{1}{2} \}} \leq \gamma |b_{i}-a_{i}| \ ; \ |
    \mathcal{R}^{t_{i}}_{a_{i},b_{i}}| > \delta |b_{i}-a_{i}|\ \biggl|
    \ a_{i} \overset{\widetilde{h_{t_{i}}}}{\leftrightarrow} b_{i}
    \biggl) 
  \leq \exp\big(-\varphi|b_{i}-a_{i}|\big). \label{L:gammalemma}
\end{equation}
\end{Lemma}

\begin{proof}
For some $\gamma>0$ to be specified later, we write
\begin{align} \label{E:Udecomp}
&P\biggl( \sum_{k=2}^{N}
  \delta_{\{|\Delta (r_{k})| \geq \frac{1}{2} \}} \leq
  \gamma |b_{i}-a_{i}| \ ; \ |
  \mathcal{R}^{t_{i}}_{a_{i},b_{i}}| > N \biggl| \
  a_{i} \overset{\widetilde{h_{t_{i}}}}{\leftrightarrow} b_{i} \biggl)\\
&\leq \sum_{(z_{1},\cdots,z_{R}) \in \mathcal{U}}
  P\bigg(\mathcal{Q}^{t_{i}}_{a_{i},b_{i}} \subset
  \bigcup_{j=1}^{R}\mathcal{H}_{z_{j}}^{t_{i}}\  ;
  \sum_{k=2}^{N}
  \delta_{\{|\Delta (r_{j})| \geq \frac{1}{2} \}} \leq
  \gamma |b_{i}-a_{i}| \ \biggl| \
  a_{i} \overset{\widetilde{h_{t_{i}}}}{\leftrightarrow} b_{i} \biggl)
  \notag \\
&\leq \sum_{(z_{1},\cdots,z_{R}) \in \mathcal{U}} P\bigg(
  \mathcal{Q}^{t_{i}}_{a_{i},b_{i}} \subset
  \bigcup_{j=1}^{R}\mathcal{H}_{z_{j}}^{t_{i}} ;
  \sum_{j=2}^{R}
  \delta_{\{|\widetilde{\Delta} (\mathcal{H}^{t_{i}}_{z_{j}})| \geq
  \frac{1}{2} \}} \leq \gamma |b_{i}-a_{i}| \ \biggl| \
  a_{i} \overset{\widetilde{h_{t_{i}}}}{\leftrightarrow} b_{i} \biggl)
  \notag \\
&\leq \sum_{(z_{1},\cdots,z_{R}) \in \mathcal{U}}
   P\bigg(\mathcal{Q}^{t_{i}}_{a_{i},b_{i}} \subset
  \bigcup_{j=1}^{R}\mathcal{H}_{z_{j}}^{t_{i}}\biggl| \
  a_{i} \overset{\widetilde{h_{t_{i}}}}{\leftrightarrow} b_{i}
  \biggl)\times \notag \\ 
& \qquad \qquad P\bigg(\sum_{j=2}^{R}
  \delta_{\{|\widetilde{\Delta} (\mathcal{H}^{t_{i}}_{z_{j}})| \geq
  \frac{1}{2} \}}
  \leq \gamma |b_{i}-a_{i}| \ \biggl| \ 
   \mathcal{Q}^{t_{i}}_{a_{i},b_{i}} \subset
  \bigcup_{j=1}^{R}\mathcal{H}_{z_{j}}^{t_{i}}\ ; \  a_{i}
  \overset{\widetilde{h_{t_{i}}}}{\leftrightarrow} b_{i}  \bigg). \notag
\end{align}

We will bound the second probability in the last sum. In order to
do this, we will describe a ``renewal shifting'' procedure. For
$\omega \in \{\mathcal{Q}^{t_{i}}_{a_{i},b_{i}} \subset
\cup_{j=1}^{R}\mathcal{H}_{z_{j}}^{t_{i}}\ ; \  a_{i}
\overset{\widetilde{h_{t_{i}}}}{\leftrightarrow} b_{i}
\}$, satisfying $|\widetilde{\Delta}
(\mathcal{H}^{t_{i}}_{z_{j}})| < \frac{1}{2} $ for some fixed
$j\geq2$, this procedure will produce a configuration
$\widetilde{\omega} \in \{\mathcal{Q}^{t_{i}}_{a_{i},b_{i}} \subset
\cup_{j=1}^{R}\mathcal{H}_{z_{j}}^{t_{i}}\ ; \  a_{i}
\overset{\widetilde{h_{t_{i}}}}{\leftrightarrow} b_{i} \}$, which
has at most a bounded number of bonds different from $\omega$, and which
satisfies
$|\widetilde{\Delta} (\mathcal{H}^{t_{i}}_{z_{j}})| \geq
\frac{1}{2} $. Moreover, this procedure maps
at most $2^{m}$ configurations to the same $\widetilde{\omega}$,
where $m$ is the number of possibly-adjusted bonds. Once this procedure
is described, for constants $c_{1},c_{2},\cdots,c_{j-1}$ we get
\begin{align} \label{E:lambdarelation}
P &\biggl( |\widetilde{\Delta} (\mathcal{H}^{t_{i}}_{z_{j}})|<1/2
  \ \biggl| \  \mathcal{Q}^{t_{i}}_{a_{i},b_{i}} \subset
  \bigcup_{j=1}^{R}\mathcal{H}_{z_{j}}^{t_{i}}\ ; \  a_{i}
  \overset{\widetilde{h_{t_{i}}}}{\leftrightarrow} b_{i}\  ;  \
  \widetilde{\Delta} (\mathcal{H}^{t_{i}}_{z_{k}})=c_{k}, \
  1\leq k < j
  \biggl) \\
&\leq \lambda'  \  P\biggl(|\widetilde{\Delta}
  (\mathcal{H}^{t_{i}}_{z_{j}})|\geq 1/2 \ \biggl| \ 
  \mathcal{Q}^{t_{i}}_{a_{i},b_{i}} \subset
  \bigcup_{j=1}^{R}\mathcal{H}_{z_{j}}^{t_{i}}\ ; \  a_{i}
  \overset{\widetilde{h_{t_{i}}}}{\leftrightarrow} b_{i} \ ; \ 
  \widetilde{\Delta} (\mathcal{H}^{t_{i}}_{z_{k}})=c_{k}, \ 1\leq k
  < j \biggl), \notag
\end{align}
where $\lambda'=\lambda'(p)>0$. This yields
\begin{align}
P\biggl(&|\widetilde{\Delta} (\mathcal{H}^{t_{i}}_{z_{j}})|\geq 1/2
  \ \biggl| \ \mathcal{Q}^{t_{i}}_{a_{i},b_{i}} \subset
  \bigcup_{j=1}^{R}\mathcal{H}_{z_{j}}^{t_{i}}\ ; \  a_{i}
  \overset{\widetilde{h_{t_{i}}}}{\leftrightarrow} b_{i} ; \ 
  \widetilde{\Delta} (\mathcal{H}^{t_{i}}_{z_{k}})=c_{k}, \mbox{ for
  } 1\leq k <j \biggl) \notag \\ 
&\geq  \frac{1}{1+ \lambda'} \notag
\end{align}
which is sufficient to bound the last probability in (\ref{E:Udecomp})
by
$P\bigl(X<\gamma |b_{i}-a_{i}|  ) ,$
where $X$ is binomially distributed with parameters 
$R-1$ and
$p^{*}=\frac{1}{1+\lambda'}$.
Taking $\gamma < p^*$ and using a bound from
\cite{Ho} we have
\begin{eqnarray*}
P\bigl(X<\gamma |b_{i}-a_{i}|\big) \leq  \exp \biggl(-
\frac{(R-1)(p^{*} - \gamma)^{2}}{2}\bigg)
 \leq
\exp(-\varphi|b_{i}-a_{i}| \bigl),
\end{eqnarray*}
for some $\varphi >0$. Using this in the right side of (\ref{E:Udecomp})
and observing that the events $\{\mathcal{Q}^{t_{i}}_{a_{i},b_{i}}
\subset
\bigcup_{j=1}^{R}\mathcal{H}_{z_{j}}^{t_{i}} \}$ are disjoint for
distinct $(z_{1},z_{2},\cdots,z_{R})\in \mathcal{U}$, 
we obtain$~\eqref{L:gammalemma}$, after summing over all
$(z_{1},z_{2},\cdots,z_{R})\in \mathcal{U}$ .

 \par
The proof will be completed by description of the ``renewal shifting''
procedure. For a given configuration $\omega\in 
\{\mathcal{Q}^{t_{i}}_{a_{i},b_{i}} \subset
\bigcup_{j=1}^{R}\mathcal{H}_{z_{j}}^{t_{i}}\ ; \  a_{i}
\overset{\widetilde{h_{t_{i}}}}{\leftrightarrow} b_{i}\}$ and a fixed
$j \leq R$, let us assume $|\widetilde{\Delta}
(\mathcal{H}^{t_{i}}_{z_{j}})|<\frac{1}{2}$, for some $j$. We will
define $\widetilde{\omega}$ by modifying some dual bonds inside
$\mathcal{S}_{z_{j}-4e_i,z_{j}+4e_i}^{t_{i}}$. Since $t_{i}$ has
slope $\frac{r}{q}$, there exists infinitely many equally spaced lattice
points on the line
$\mathcal{H}_{z_{j}}^{t_{i}}$. We will use  one of the two lattice
points on $\mathcal{H}_{z_{j}}^{t_{i}}$ closest to the regeneration
point $\sigma_{j}$.  Call these locations $u_{j}$ and
$v_{j}$, with $u_{j}=\sigma_{j} + (-r,q)$ and
$v_{j}=\sigma_{j}+(r,-q)$. The configuration $\omega$ has open dual bonds
$\langle \sigma_{j}-e_i,\sigma_{j}\rangle$ and $\langle
\sigma_{j},\sigma_{j}+e_i\rangle$. 

There exists a path $\gamma^L_j$ from $\sigma_{j}-2e_i$ to $u_j - e_i$
in $\mathcal{S}^{t_{i}}_{z_{j}-3e_i,z_{j}-e_i}$ having all steps upward
or leftward, with $\gamma^L_j \cap \mathcal{H}^{t_i}_{\sigma_j - e_i}
= \{u_j - e_i\}$, and similarly a path 
$\gamma^R_j$ from $\sigma_{j}+2e_i$ to $u_j + e_i$
in $\mathcal{S}^{t_{i}}_{z_{j}+e_i,z_{j}+3e_i}$ having all steps upward
or leftward with $\gamma^R_j \cap \mathcal{H}^{t_i}_{\sigma_j + e_i}
= \{u_j + e_i\}$.  Let $A_j$ be the closed region
bounded by $\gamma^L_j,
\gamma^R_j$ and the horizontal lines through $\sigma_j$ and $u_j$.  To
make our choice of $\gamma^L_j,
\gamma^R_j$ unique, let us specify that $A_j$ be maximal under the
constraints we have imposed on $\gamma^L_j,
\gamma^R_j$.  Let
$\mathcal{D}_j$ be the set of all dual bonds having one endpoint in
$\partial A_j$ and the other outside $A_j$.  Note there are at most $12
q$ dual bonds contained in $A_j$, and at most $2r+2q+10$ dual bonds in
$\mathcal{D}_j$.  Let
$\widetilde{\omega}$ be such that
\begin{enumerate}
\item[(1)] all dual bonds in $\partial A_j \backslash \{ 
\langle \sigma_{j}-e_i,\sigma_{j}\rangle, \langle
\sigma_{j},\sigma_{j}+e_i\rangle \}$ are open;
\item[(2)] all other dual bonds contained in $A_j$ are closed;
\item[(3)] all dual bonds in $\mathcal{D}_j \cap
\mathbf{C}_{a_i,b_i}^{t_i}(\omega)$ are open;
\item[(4)] all dual bonds in $\mathcal{D}_j \backslash
\mathbf{C}_{a_i,b_i}^{t_i}(\omega)$ are closed;
\item[(5)] all other dual bonds retain their state from $\omega$.
\end{enumerate}
In the altered configuration $\widetilde{\omega}$, the regeneration
point is still on $\mathcal{H}_{z_{j}}^{t_{i}}$ but shifted from
$\sigma_{j}$ to $u_{j}$.  After these alterations, if
$|\widetilde{\Delta} (\mathcal{H}^{t_{i}}_{z_{j}})|\geq\frac{1}{2}$,
then we are done. It is possible that $|\widetilde{\Delta}
(\mathcal{H}^{t_{i}}_{z_{j}})| < \frac{1}{2}$, for the following reason.
Let $k$ be such that $z_j = r_k$.  If there are other
regeneration points in $\mathcal{S}^{t_{i}}_{z_{j}-3e_i,z_{j}+3e_i}$ in
$\omega$, shifting the regeneration point to $u_{j}$ will destroy these
regeneration points; any regeneration points in 
$\mathcal{S}^{t_{i}}_{z_{j}-4e_i,z_{j}+4e_i} \backslash 
\mathcal{S}^{t_{i}}_{z_{j}-3e_i,z_{j}+3e_i}$ in $\omega$ may or may not
be destroyed, depending on the exact geometry of the situation.
At any rate, if $r_{k-1}$ is destroyed, the new
``preceding regeneration point'' for $z_j$ will be
outside the slab $\mathcal{S}^{t_{i}}_{z_{j}-3e_i,z_{j}+3e_i}$, equal to
$r_{k-2}$ or $r_{k-3}$, and we may have $|\widetilde{\Delta}
(\mathcal{H}^{t_{i}}_{z_{j}})|<\frac{1}{2}$ in $\widetilde{\omega}$,
depending on the location of this new preceding regeneration point
relative to
$l_i$.  If this is the case we shift the regeneration point from 
$\sigma_{j}$ to $v_{j}$ instead of
$u_{j}$.  For this we use paths $\tilde{\gamma}^L_j$ from
$\sigma_{j}-2e_i$ to $v_j - e_i$ and $\tilde{\gamma}^R_j$ from
$\sigma_{j}+2e_i$ to $v_j + e_i$ in place of $\gamma^L_j$ and
$\gamma^R_j$, under an analogous maximality constraint.  Let $x^L_j$
(respectively
$x^R_j$) be the site in
$\gamma^L_j$ (respectively $\gamma^R_j$) closest to
$\mathcal{H}^{t_i}_{\sigma_j - 3e_i}$ (respectively
$\mathcal{H}^{t_i}_{\sigma_j + 3e_i}$).  Due to the maximality
constraints we have imposed, since all our slabs have boundaries with
slope $-q/r$, 
$x^L_j +(r,-q)$ is the site in $\tilde{\gamma}^L_j$ closest to 
$\mathcal{H}^{t_i}_{\sigma_j - 3e_i}$, and $x^R_j +(r,-q)$ is the site in
$\tilde{\gamma}^R_j$ closest to 
$\mathcal{H}^{t_i}_{\sigma_j + 3e_i}$.  This means that
$\gamma^L_j$ and $\tilde{\gamma}^L_j$ intersect the same slabs
orthogonal to $t_i$, and similarly for $\gamma^R_j$ and
$\tilde{\gamma}^R_j$.  As a consequence, the same regeneration points
are destroyed, regardless of whether we shift to $u_j$ or $v_j$, so
$\widetilde{\omega}$ has the same preceding regeneration point either
way.  It follows that if shifting to $u_j$ results in 
$|\widetilde{\Delta}
(\mathcal{H}^{t_{i}}_{z_{j}})|<\frac{1}{2}$, then shifting to $v_j$
results in $|\widetilde{\Delta}
(\mathcal{H}^{t_{i}}_{z_{j}})| \geq \frac{1}{2}$, i.e. there is always a
shift (the one we choose to create $\widetilde{\omega}$) which results
in $|\widetilde{\Delta}
(\mathcal{H}^{t_{i}}_{z_{j}})| \geq \frac{1}{2}$.

Note that only a bounded
number of different configurations may map to the same configuration
$\widetilde{\omega}$. In any case, $\widetilde{\omega}$ and $\omega$
yield the same value of $\mathcal{Q}^{t_{i}}_{a_{i},b_{i}}$, and the
probabilities of $\omega$ and
$\widetilde{\omega}$ are within a bounded factor (depending on $p$),
which yields$~\eqref{E:lambdarelation}$, completing the proof.
\end{proof}

\section{Exchangeability of Increments}
The core idea in our proof of $\eqref{E: mainthm}$ is to make
use of the renewal structure of connectivities, for connections
between any two consecutive extreme points $w_i,w_{i+1}$ in the
$s$-hull skeleton  with $i \in
\mathcal{L}$, to see that the increments $\Delta(r_j), 2 \leq j \leq N$,
form an exchangeable sequence under certain conditioning, that is, the
joint distribution is permutation invariant.  The partial sums of this
sequence behave like those of an i.i.d. sequence, and from this we can
show that with high probability, the path of open dual bonds will not
stay in the ``narrow tube'' from $w_i$ to
$w_{i+1}$ with diameter
$2d=4\theta l^{1/3}(\log l)^{-2/3}$.  
In this section we will prove
this exchangeability.  Let 
$\epsilon, t_i, a_i, b_i$ be as in the preceeding, $\delta$ as in
$\eqref{L: noofrenewals}$ and $\gamma$ as in
$\eqref{L:gammalemma}$. Define the event
$E=E(a_{i},b_{i},t_{i},\gamma,\delta)$ by
\begin{eqnarray*}
   E=\bigg\{ a_{i}
  \overset{\widetilde{h_{t_{i}}}}{\longleftrightarrow} b_{i}\bigg\} \cap
  \biggl\{
  |\mathcal{R}^{t_{i}}_{a_{i},b_{i}}|\geq\delta|b_{i}-a_{i}|\biggr\}
  \cap
  \biggl\{ \sum_{k=2}^{N}
  \delta_{\{|\Delta(r_{k})| \geq \frac{1}{2} \}} \geq
  \gamma |b_{i}-a_{i}|\biggr\}.
\end{eqnarray*}
For $v,w \in
T_{d}^{i} \cap \mathcal{S}_{a_{i},b_{i}}^{t_{i}} \cap
(\mathbb{Z}^{2})^*$, define the sets
\begin{eqnarray*}
  V(v,w)= \biggl\{ \ \zeta=(\zeta_{2},...,\zeta_{N})\in
  \mathbb{R}^{N-1}:
  |\zeta_{2}|\geq |\zeta_{3}| \geq \cdots \geq |\zeta_{N}| \ ; 
  \sum_{k=2}^{N}
  \zeta_{i}= f(w)-f(v) \biggl\}.
\end{eqnarray*}
For given $\zeta' \in V(v,w)$, let $ F = F(v,w,\zeta')$ denote the event
that the following all hold:
\begin{enumerate} 
\item[(i)] $ a_{i}
\overset{\widetilde{h_{t_{i}}}}{\longleftrightarrow} b_{i}$ ,
\item[(ii)] the first  and $N$-th
regeneration points are at $v$ and $w$, respectively,
\item[(iii)] for some permutation $\pi:\{2,\cdots,N\} \to
\{2,\cdots,N\}$, we have  
\[
  \Delta (r_{2})=\zeta'_{\pi(2)},
  \Delta (r_{3})=\zeta'_{\pi(3)}, \cdots, \Delta (r_{N})=\zeta'_{\pi(N)}.
\]
\end{enumerate}
Observe that condition (iii) determines the values of the $\Delta
(r_{k})$'s up to an ordering, and (ii) and (iii) imply $\sum_{k=2}^{N}
\Delta(r_{k})= f(w)-f(v)$.

\begin{Lemma} 
For fixed $a_{i},b_{i},t_{i},\gamma,\delta,v,w,\zeta',E,F$ as in the
preceding,
$\Delta(r_{2}),\cdots, \Delta(r_{N})$ are
exchangeable under the measure $P(\cdot \ | \ E \cap F)$.
\end{Lemma}

\begin{proof}
$E\cap F $ determines the location of 
first and $N$-th regeneration points,
and values of increments in between them, up to an ordering. We
will first show how to exchange any two adjacent increments.
Consider a configuration $\omega \in E\cap F$, with $
\Delta(r_{2})=\zeta'_{2},\Delta(r_{3})=\zeta'_{3},\cdots ,
\Delta(r_{N})=\zeta'_{N}$. For fixed $k\geq 2$, let us consider
increments $\Delta(r_{k})$ and $\Delta(r_{k+1})$.  By definition of regeneration points,
the bonds that are only partially in the slab
$\mathcal{S}^{t_i}_{r_k,r_{k+1}}$ or have exactly one endpoint in the
within-slab cluster containing $r_k$ and $r_{k+1}$ are all vacant. We
construct a configuration $\widetilde{\omega}$ such that outside
$\mathcal{S}^{t_{i}}_{r_{k-1},r_{k+1}}$ we have
$\widetilde{\omega}=\omega$. We obtain $\widetilde{\omega}$ by
interchanging the relative positions of the configurations
$\omega_{\mathcal{S}^{t_{i}}_{r_{k-1},r_{k}}}$ and
$\omega_{\mathcal{S}^{t_{i}}_{r_{k},r_{k+1}}}$ and moving the
bonds crossing $\mathcal{H}_{r_{k}}^ {t_{i}}$ so that they cross
$\mathcal{H}_{r_{k-1}+(r_{k+1}-r_{k})}^ {t_{i}}$ instead.  The latter
move is done in such a way that the relative positions of the bonds
remain the same, with the old position relative to $r_k$ becoming the
new position relative to $r_{k-1} + (r_{k+1} - r_k)$.  More precisely,
the configuration in $\mathcal{S}^{t_{i}}_{r_{k-1},r_{k}}$
is translated by $r_k - r_{k-1}$, the configuration in 
$\mathcal{S}^{t_{i}}_{r_{k},r_{k+1}}$ is translated by $r_k -
r_{k+1}$, and
each bond touching or crossing $\mathcal{H}_{r_{k}}^ {t_{i}}$ is
translated by
$r_{k-1} + (r_{k+1} - 2r_k)$.  This moves the
$k$-th regeneration point from $r_k$ to $r_{k-1} + (r_{k+1} - r_k)$,
without altering the locations of other regeneration points.  
The configuration $\widetilde{\omega}$ is in $E\cap F$ and the
increments of $\widetilde{\omega}$ satisfy
\begin{eqnarray*}
\Delta(r_{k})=\zeta'_{j+1},
\Delta(r_{k+1})=\zeta'_{j},\mbox{ and }\Delta(r_{m})=\zeta'_{m}, \mbox{
for } 2\leq m
\leq N,\ k \ne m,m+1.
\end{eqnarray*}
Moreover, replacing $\omega$ with $\widetilde{\omega}$
does not affect probability under the measure $P(\cdot|\ E\cap F)$,
due to shift invariance.  We can repeat the exchanging
of adjacent increments until we achieve the desired permutation of
$\{2,3\cdots,N\}$, and the lemma follows.
\end{proof}

\section{Staying in the Narrow Tube}
In this section, we will show that there is an extra probabilistic
cost associated to the event that $\Gamma_{0}$ stays in the
narrow tube $T_{d}^{i}$, between $a_{i}$ and $b_{i}$ . The proof
involves randomization of the order of the increments, using
exchangeability.  

\begin{Lemma} 
Let $i\in \mathcal{L}$ and let $\epsilon, t_i, a_i, b_i$
be as in the preceeding.  Let $\delta$ be as in (\ref{L: noofrenewals})
and $\gamma$ as in (\ref{L:gammalemma}). There exists
$\kappa=\kappa(\gamma)>0$  such that
for all $v,w\in T_{d}^{i}\cap
\mathcal{S}^{t_{i}}_{a_{i},b_{i}} \cap
(\mathbb{Z}^{2})^*$ and $\zeta' \in V(v,w) \cap [-2d,2d]^{N}$, for
$E=E(a_{i},b_{i},t_{i},\gamma,\delta)$, $F=F(v,w,\zeta')$, provided $l$
is large we have
\begin{eqnarray} \label{E:intmainlemma}
  P\biggl(a_{i}
  \overset{\widetilde{h_{t_{i}}}}{\longleftrightarrow} b_{i} \mbox{ in }
  T_{d}^{i} \ \bigg| \ E \cap F\bigg) \leq
  2 \exp\biggl(\frac{-\kappa|w_{i+1}-w_{i}|}{d^{2}}\biggl).
\end{eqnarray}
\end{Lemma}

\begin{proof}
Observe that when $a_{i}
\overset{\widetilde{h_{t_{i}}}}{\longleftrightarrow} b_{i}$ in $
T_{d}^{i}$, every open path from $a_i$ to $b_i$ in $T_{d}^{i}$ must pass
through all regeneration points.  Thus
\begin{equation} \label{E:regenloc}
  P\big(a_{i} \overset{\widetilde{h_{t_{i}}}}{\longleftrightarrow}
  b_{i} \mbox{ in } T_{d}^{i} \ | \ \ E\cap F \big)  \leq
  P(\{r_{1},r_{2},\cdots,r_{N}\} \subset T_{d}^{i} \cap
  \mathcal{S}^{t_{i}}_{a_{i},b_{i}} \ \big| \ E\cap F \big).
\end{equation}
We can relate the last probability to an event
involving increments. If $1\leq k_{1}<
k_{2}\leq N$ and $\big|\sum_{j=k_{1}}^{k_{2}-1}
\Delta(r_{j+1})\big| > 2d$, then the $k_{1}$-th or
$k_{2}$-th regeneration point must lie outside of $T_{d}^{i}$.
Therefore,
\begin{eqnarray}
& &P(\{r_{1},r_{2},\cdots,r_{N+1}\}
  \in T_{d}^{i}
  \cap \mathcal{S}^{t_{i}}_{a_{i},b_{i}} \ \big| \ E\cap F \big)
  \nonumber \\ \label{E:blocks1}
&\leq& P\bigg( \bigg|
  \sum_{j=k_{1}}^{k_{2}-1} \Delta(r_{j+1})\bigg| \leq 2d, \mbox{ for all
  } 1\leq k_{1} < k_{2} \leq N  \ \bigg|  \ E\cap F \bigg).
\end{eqnarray}
Instead of looking at 
partial sums for all possible values of $k_{1},k_{2}$, we will consider
disjoint blocks of increments with random lengths
$X_{1},X_{2},\cdots, X_{B}$ satisfying
$X_{1}+\cdots+X_{B}< N$, for some $B\in
\mathbb{N}$. Let $S_{n}=\sum_{k=1}^{n}X_{k}$, for $n=1,\cdots ,B$,
and let $S_{0}=0$. Then$~\eqref{E:blocks1}$ is bounded by
\begin{eqnarray} \label{E:blocks2}
P\bigg(\bigcap_{k=1}^{B}\bigg\{\max_{1\leq m \leq
X_{k}}\bigg|\sum_{j=S_{k-1}+1}^{S_{k-1}+m} \Delta(r_{j+1})\bigg|
\leq 2d\bigg\}  \bigg|   \ E\cap F
\bigg),
\end{eqnarray}
where we define $X_{0}=0$.  If we take $X_{k}$, for  $1\leq k \leq
B$, to be deterministic, the increments on
these disjoint blocks will not be independent of the increments
on other blocks. In order to
reduce the dependence between these disjoint blocks, we will
take the $X_{k}$'s to be (non-independent) binomially distributed random
variables. Next, we use
exchangeability of increments to write$~\eqref{E:blocks2}$ in an
equivalent form. For binomially distributed $X_{1}$ with
parameters, $N-1$ and $p_{0}$, with $p_{0}$ to be specified
later,
\begin{eqnarray*}
  \sum_{j=2}^{X_{1}+1}\Delta(r_{j})  \overset{d}{=}
  \sum_{j=2}^{N}\delta_{j1} \zeta'_{j}
\end{eqnarray*}
where the $\delta_{j1},\ j = 2,..,N$ are i.i.d. Bernoulli random
variables with parameter $p_{0}.$ That is, the sum of
first $X_{1}$ increments have the same distribution as the sum of
increments randomly selected according to the $\delta_{j1}$'s.
Continuing this way, for each following random block, we
replace the sum of increments corresponding to that block with a sum
of increments that are randomly selected from those increments
remaining after the earlier steps of the increment--selection process.
More precisely, we do the following.  Define $p_{0}$ and the number of
blocks by
\begin{equation*}
  p_{0} = \frac{K_{19}d^{2}}{|w_{i+1}-w_{i}|},
  \qquad B = \lfloor \frac{1}{2p_{0}} \rfloor,
\end{equation*}
 where $K_{19}=K_{19}(\gamma)$ is sufficiently large constant, to
be specified later. Observe that $p_0 = O((\log
l)^{-1})$.  For all $2 \leq j \leq N$, define
\begin{eqnarray} \label{E:deltas}
  \delta_{jk}=\begin{cases} 0 & \mbox{ with probability }
  p_{k}=\frac{1-kp_{0}}{1-(k-1)p_{0}} \\ 1 & \mbox{ with probability
  } 1-p_{k}=\frac{p_{0}}{1-(k-1)p_{0}}, \end{cases}
\end{eqnarray}
with $\{\delta_{jk}, \ j=2,\cdots N\ $, $k=1,\cdots,B\}$ 
independent random variables. Also define
$Y_{jk}=(1-\delta_{j1})(1-\delta_{j2})\cdots(1-\delta_{jk})$, for
$j=2,..,N$,  $k=1,..,B$.  Then we have
\begin{eqnarray} \label{E:Yjks}
  Y_{jk}=\begin{cases} 0 & \mbox{ with probability } kp_{0} \\
  1 & \mbox{ with probability } 1-kp_{0} \end{cases}
\end{eqnarray}
The random variable $Y_{jk}=1$ says that the $j$-th increment is
not selected for the first $k$ blocks, and
$Y_{j(k-1)}\delta_{jk} = 1$ says that the $j$-th
increment is selected for the $k$-th block. We define the
length of the $k$-th block $X_{k}$ as
\begin{eqnarray*} 
  X_{k}=\sum_{j=2}^{N} Y_{j(k-1)}\delta_{jk}, \quad \mbox{ for
  }k=1,2,\cdots,B.
\end{eqnarray*}
It can be easily
seen from$~\eqref{E:deltas}$ and$~\eqref{E:Yjks}$
that the
$X_{k}$'s are binomially distributed (but not independent) with
parameters $N-1$ and $p_{0}$. 
By exchangeability, we can
rewrite$~\eqref{E:blocks2}$ as
\begin{eqnarray} \label{E:blocks2b}
  P\bigg(\bigcap_{k=1}^{B}\bigg\{\max_{2\leq m \leq N}
  \bigg|\sum_{j=2}^{m}Y_{j(k-1)}\delta_{jk}\zeta_{j}\bigg|\leq 2d
  \bigg\} \ \bigg| \ E \cap F \bigg),
\end{eqnarray}
since $\sum_{k=1}^{B}X_k \leq N$, which holds because no $j$ can
be chosen for more than one block.  We let
\begin{eqnarray*}
  D_{k} = \bigg\{ \omega :\max_{2\leq m \leq N}
  \bigg| \sum_{j=2}^{m}Y_{j(k-1)}\delta_{jk}\zeta_{j}\bigg|\leq 2d
  \bigg\}.
\end{eqnarray*}
We need to control the number of
increments $|\zeta_j| \geq 1/2$ 
which remain after some blocks have been selected. By definition of
$E$, there are at least $\lfloor\gamma|b_{i}-a_{i}|\rfloor$ such
increments before the first block
is selected. Let $g=\sqrt{B\lfloor \gamma|b_{i}-a_{i}|\rfloor}$, let
$G_{1}=E\cap F$, and for $k=2,\cdots,B$ define
\begin{eqnarray*}
   G_{k}= \bigg\{\omega:
  \bigg|\bigg(\sum_{j=2}^{\lfloor
  \gamma|b_{i}-a_{i}|\rfloor}Y_{j(k-1)}\bigg)-(1-(k-1)p_{0})\lfloor
  \gamma|b_{i}-a_{i}|\rfloor\bigg|\leq g \bigg\}.
\end{eqnarray*}
Let $I_{k}=\{j: Y_{j(k-1)}=1, 1\leq j \leq N-1 \}$, be the random
set of remaining increment indices before the $k$-th block is
selected. Then $G_{k-1}$ provides control over $\big|
I_{k}\cap \{1,2,\cdots,\lfloor \gamma|b_{i}-a_{i}|\rfloor\}
\big|$, the number of remaining increments that are greater than
or equal to $1/2$; here we use
the monotonicity of the $|\zeta_{j}|$'s, and the fact that at least
$\lfloor
\gamma|b_{i}-a_{i}|\rfloor$ increments are greater than or equal
to $1/2$.
 By$~\eqref{E:regenloc}$--$\eqref{E:blocks2}$ we have
\begin{eqnarray}
& & P\big(a_{i}
  \overset{\widetilde{h_{t_{i}}}}{\longleftrightarrow} b_{i} \mbox{
  in } T_{d}^{i}\ \big| \ E\cap F \big)  \nonumber \\ \label{E:blocks3}
& &\qquad \leq P\bigg(\bigcap_{k=1}^{B} \big(D_{k}\cap G_{k}
  \big) \ \bigg| \ E
  \cap F\bigg) + P\bigg( \bigg[\bigcap_{k=1}^{B} G_{k}\bigg]^{c} \
  \bigg| \ E \cap F \bigg)
\end{eqnarray}
First we bound the probability in (\ref{E:blocks3}) of a large deviation
for some block for the number of available large increments, using a
bound from \cite{Ho}:
\begin{align}
  P\bigg( \bigg[\bigcap_{k=1}^{B} G_{k}\bigg]^{c} \ \bigg| \ E \cap
    F \bigg) &\leq \sum_{k=2}^{B}P(G_{k}^{c} \mid E \cap F) \\
  &\leq 2B \exp\big(\frac{-2g^{2}}{\lfloor
    \gamma|b_{i}-a_{i}|\rfloor}\big) \notag \\
  &= 2B \exp(-2B) \notag \\
  &\leq \exp(-B), \notag
\end{align}
where the last inequality holds for $l$ sufficiently large.
Next, we consider the probability the probability of staying in the
narrow tube in the absence of such a large deviation.  This probability
from (\ref{E:blocks3}) can be written
\begin{align}
  P\big( &D_{1} \ \big|\  E \cap F\big) \times \prod_{k=2}^{B}
    P\bigg(D_{k}\cap G_{k}  \ \bigg| \ E \cap F \cap
    \bigcap_{j=1}^{k-1}(D_{j}\cap G_{j} )\bigg) \label{E:lastprod} \\
  &\leq P\big(D_{1} \ \big|\  E \cap F\big) \times \prod_{k=2}^{B}
    P\bigg(D_{k} \ \bigg| \ E \cap F \cap
    \bigcap_{j=1}^{k-1}(D_{j}\cap G_{j})\bigg). \notag
\end{align}
We will conclude by showing
\begin{eqnarray} \label{E:laststep}
  P\bigg(D_{k} \ \bigg| \ E \cap F \cap
  \bigcap_{j=1}^{k-1}(D_{j}\cap G_{j} )\bigg) \leq 2/3,
\end{eqnarray}
for $k\geq 2$. The proof that  $\  P\big(D_{1} \ \big|\  E \cap
F\big) \leq 2/3 \  $ follows by the same technique. Let us fix $k\geq
2$.  We define a family of
sets of indices:
\begin{equation*}
  \mathcal{I}_{k}= \bigg\{\Upsilon\subset \{2,\cdots,N-1\} : \bigg|
  \sum_{j=2}^{\lfloor \gamma|b_{i}-a_{i}|\rfloor} \delta_{\{j \in
  \Upsilon\}} - (1-(k-1)p_{0})\lfloor
  \gamma|b_{i}-a_{i}|\rfloor\bigg| \leq g \bigg\}.
\end{equation*}
 For $\Upsilon \in \mathcal{I}_{k}$, and $n\leq
N-1$, define
\begin{eqnarray*}
  \Upsilon_{n}=\Upsilon \cap \{2,3, \cdots, n\}
\end{eqnarray*}
Observe that if $G_{k-1}$ occurs then $I_{k} \in
\mathcal{I}_{k}$.  It follows that
\begin{eqnarray} \label{E:step1}
& &P\bigg(D_{k} \ \bigg| \ E \cap F \cap
  \bigcap_{j=1}^{k-1}(D_{j}\cap G_{j} )\bigg) \nonumber \\ 
& & \qquad = \sum_{\Upsilon \in \mathcal{I}_{k}} P\bigg(D_{k} \cap
  \{ I_{k}=\Upsilon \} \ \bigg| \ E \cap F \cap
  \bigcap_{j=1}^{k-1}(D_{j}\cap G_{j})\bigg).
\end{eqnarray}
Fix $\Upsilon \in \mathcal{I}_k$ and define the event $H_{k}=[I_{k} =
\Upsilon ]  \cap E \cap F
\cap \bigcap_{j=1}^{k-1}(D_{j}\cap G_{j})$.  Define
\begin{eqnarray*}
  Q(k,\Upsilon_{n})=
  \Bigg[\frac{1}{2}\Bigg(\frac{p_{0}(1-kp_{0})}{(1-(k-1)p_{0})^{2}}
  \Bigg)
  \sum_{j\in \Upsilon_{n}}(\zeta'_{j})^{2}\Bigg]^{1/2},
\end{eqnarray*}
so that
\begin{eqnarray*}
  \mbox{Var}\bigg(\sum_{j\in
  \Upsilon_{n}}\delta_{jk}\zeta'_{j}\ \bigg|\ H_k \bigg)
  =2[Q(k,\Upsilon_{n})]^{2}.
\end{eqnarray*}
\par
 For any index set
$\Upsilon \in \mathcal{I}_{k}$,  one of three possibilities has
to hold:
\begin{enumerate}
\item[(1)] for all $n=2,3,...,N$
\begin{eqnarray*}
  \bigg|\mathbb{E}\bigg(\sum_{j\in
  \Upsilon_{n}}\delta_{jk}\zeta'_{j}\ \bigg|\ H_k \bigg)\bigg|\ = \
  \ \bigg| \sum_{j\in \Upsilon_{n}}
  \frac{p_{0}}{1-(k-1)p_{0}}\zeta'_{j}
  \bigg|&\leq&  2d+ Q(k,\Upsilon_{n});
\end{eqnarray*}
\item[(2)] for some $n_{0}$, $2\leq n_{0}\leq N$
\begin{eqnarray*}
  \mathbb{E}\bigg(\sum_{j\in
  \Upsilon_{n_{0}}}\delta_{jk}\zeta'_{j}\ \bigg|\ H_k \bigg) &>&
  2d+Q(k,\Upsilon_{n_{0}});
\end{eqnarray*}
\item[(3)] for some $n_{0}$, $2\leq n_{0}\leq N$
\begin{eqnarray*}
  -\mathbb{E}\bigg(\sum_{j\in
  \Upsilon_{n_{0}}}\delta_{jk}\zeta'_{j}\ \bigg|\ H_k \bigg) &>&
  2d+Q(k,\Upsilon_{n_{0}}).
\end{eqnarray*}
\end{enumerate}
In case (2), 
\begin{eqnarray*}
& &P\big(D_{k} \ \big| \ H_k \big) \\
&  &\qquad \leq P\bigg(-2d \leq \sum_{j\in \Upsilon_{n_{0}}}
  \delta_{jk}\zeta'_{j} \leq 2d \ \bigg|\ H_k \bigg)  \\ 
&  &\qquad \leq P\bigg(\sum_{j\in \Upsilon_{n_{0}}} \bigg[
  \delta_{jk}\zeta'_{j}-\frac{p_{0}}{1-(k-1)p_{0}}\zeta'_{j} \bigg] <
  -Q(k,\Upsilon_{n_{0}}) \bigg| \ H_k \  \bigg)
\end{eqnarray*}
By Chebyshev's
inequality, the last probability is bounded by
\begin{eqnarray*}
\frac{1}{1+\frac{[Q(k,\Upsilon_{n_{0}})]^{2}}
  {2[Q(k,\Upsilon_{n_{0}})]^{2}}}
  =\frac{2}{3}.
\end{eqnarray*}
In case (3), similarly, $P(D_k \mid H_k) \leq 2/3$.
Case (1) requires some extra work. Using Kolmogorov's inequality we get
\begin{align*}
P(D_{k} \mid H_k ) &\leq P\bigg(\max_{2\leq m \leq N} \Big|\sum_{j\in
  \Upsilon_{m}}
  \delta_{jk}\zeta'_{j}-\frac{p_{0}}{1-(k-1)p_{0}}\zeta'_{j}
  \Big|  < 4d+Q(k,\Upsilon_{N}) \bigg| \ H_k \  \bigg) \\
&\leq \frac{\big[6d+
  Q(k,\Upsilon_{N})\big]^{2}}{2[Q(k,\Upsilon_{N})]^{2}}.
\end{align*}
The proof of (\ref{E:laststep}) will be concluded by showing
\begin{equation*} \label{E:98bound}
  \frac{d^{2}}{2[Q(k,\Upsilon_{N})]^{2}} \leq \frac{1}{98},
\end{equation*}
since this implies
\begin{eqnarray*} 
  \frac{\big[6d+
  Q(k,\Upsilon_{N})\big]^{2}}{2[Q(k,\Upsilon_{N})]^{2}} \leq 2/3.
\end{eqnarray*}
Since $\Upsilon \in \mathcal{I}_k$ we have
\begin{eqnarray*}
\sum_{j \in \Upsilon}(\zeta'_{j})^{2} &\geq&
  \frac{1}{4}\sum_{j
  \in \Upsilon} \delta_{\big\{|\zeta'_{j}|\geq 1/2\big\}}\\
&\geq&
  \frac{1}{4}|\Upsilon_{\lfloor\gamma|b_{i}-a_{i}|\rfloor}|\\
&\geq&
  \frac{1}{4} \bigg((1-(k-1)p_{0})
  \lfloor\gamma|b_{i}-a_{i}|\rfloor-g\bigg) \\ 
&\geq& \frac{1}{8}
  \bigg( \lfloor\gamma|b_{i}-a_{i}|\rfloor-2 g\bigg) \\ 
&\geq&
  \frac{1}{16}  \gamma|b_{i}-a_{i}|
\end{eqnarray*}
since $\frac{1}{2} \leq 1-kp_{0} \leq 1$,  for sufficiently large
$l$. Therefore,
\begin{eqnarray*}
  \frac{d^{2}}{2[Q(k,N)]^{2}} \leq \frac{16d^{2}((1-(k-1)
  p_{0})^{2})}{p_{0}(1-kp_{0}) \gamma|b_{i}-a_{i}| } \
  \leq \frac{32 d^{2}}{p_{0} \gamma|b_{i}-a_{i}| } \ =
  \frac{32 |w_{i+1}-w_{i}|}{K_{19}
  \gamma|b_{i}-a_{i}|}
\end{eqnarray*}
By$~\eqref{E:wbarelation}$, we can choose
$K_{19}=K_{19}(\gamma)$ (from the definition of $p_{0}$)
sufficiently large so that the last expression is less than 1/98,
for large $l$. Under each case (1)-(3) we have shown $P(D_k \mid H_k)
\leq 2/3$, for arbitrary $\Upsilon \in \mathcal{I}_k$.
With \eqref{E:step1} this proves$~\eqref{E:laststep}$.
Using \eqref{E:blocks3}--\eqref{E:lastprod} we get
\begin{eqnarray*}
P\big(a_{i}
  \overset{\widetilde{h_{t_{i}}}}{\longleftrightarrow} b_{i} \mbox{ in }
  T_{d}^{i}\ \big| \ E\cap F \big) &\leq& (2/3)^{B} + e^{-B}\\
&\leq& 2 \exp{\big(\frac{-\kappa|w_{i+1}-w_{i}|}{d^{2}}\big)},
\end{eqnarray*}
for some $\kappa>0$, which concludes the proof of the lemma.
\end{proof}

\begin{Lemma}
Let $\epsilon, t_i, a_i, b_i$
be as in the preceeding, with $i\in \mathcal{L}$.  Let $\delta$ be as in
(\ref{L: noofrenewals}), $\gamma$ as in (\ref{L:gammalemma}) and
$\kappa$ as in (\ref{E:intmainlemma}).  Provided $l$ is sufficiently
large we have
\begin{eqnarray} \label{E:mainlemma}
  \ \ \ \ P\biggl(a_{i}
  \overset{\widetilde{h_{t_{i}}}}{\longleftrightarrow} b_{i} \mbox{
  in } T_{d}^{i} \ \bigl| \ a_{i}
  \overset{\widetilde{h_{t_{i}}}}{\longleftrightarrow} b_{i}\biggl)
  \leq  3 \exp{\biggl(\frac{-\kappa |w_{i+1}-w_{i}|}{d^{2}} \biggl)}.
\end{eqnarray}
\end{Lemma}

\begin{proof}
Let $\nu'$ be as in (\ref{L: noofrenewals2}) and $\varphi$ as in
(\ref{L:gammalemma}). We will consider intersections of the event
$\{a_{i}
\overset{\widetilde{h_{t_{i}}}}{\longleftrightarrow} b_{i} \mbox{ in }
T_{d}^{i}\}$ with the events $E=E(a_{i},b_{i},t_{i},\delta,\gamma)$
and
$E^{c}$, separately. First we have
\begin{eqnarray}
&P\biggl(E^{c} \ | \ a_{i}&
  \overset{\widetilde{h_{t_{i}}}}{\longleftrightarrow} b_{i}\biggl)
  \nonumber \\
& &\leq P\biggl(|\mathcal{R}^{t_i}_{a_{i},b_{i}}|< \delta |b_{i}-a_{i}| \
  \bigl| \
  a_{i} \overset{\widetilde{h_{t_{i}}}}{\longleftrightarrow}
  b_{i}\biggl) \ \ \ \ \ \ \ \  \ \ \ \ \ \ \ \ \ \ \ \ \ \ \ \ \ \
  \ \ \ \ \nonumber \\ 
& & \qquad +P\biggl(
  \sum_{j=1}^{N} \delta_{\{|\Delta (r_{j})| \geq \frac{1}{2} \}}
  \leq \gamma |b_{i}-a_{i}| \ ; \ |
  \mathcal{R}^{t_{i}}_{a_i,b_i}| > \delta |b_{i}-a_{i}|\ \biggl|
  \ a_{i} \overset{\widetilde{h_{t_{i}}}}{\leftrightarrow} b_{i}
  \biggl) \nonumber \\   
& & \leq \exp(-\nu'|b_{i}-a_{i}|)
  +\exp(-\varphi|b_{i}-a_{i}|) \label{E:Ec},
\end{eqnarray}
by$~\eqref{L: noofrenewals2}$ and$~\eqref{L:gammalemma}$.
Since $i\in\mathcal{L}$, this bound is small compared to the
right side of$~\eqref{E:mainlemma}$. Next, we have
\begin{eqnarray} \label{E:vwsum}
& &P\biggl(\{a_{i}
  \overset{\widetilde{h_{t_{i}}}}{\longleftrightarrow} b_{i} \mbox{ in }
  T_{d}^{i}\} \cap E
  \ \bigl| \ a_{i}
  \overset{\widetilde{h_{t_{i}}}}{\longleftrightarrow} b_{i}\biggl) \\
&\leq& \sum_{v,w\in T_{d}^{i}\cap \mathcal{S}^{t_{i}}_{a_{i}+3e_i,b_{i}}
  \cap
  (\mathbb{Z}^{2})^*} \bigg[ \sum_{\zeta' \in V(v,w)} P\biggl(\{a_{i}
  \overset{\widetilde{h_{t_{i}}}}{\leftrightarrow} b_{i} \mbox{ in }
  T_{d}^{i}\} \cap E\cap F(v,w,\zeta')
  \ \bigl| \   a_{i}
  \overset{\widetilde{h_{t_{i}}}}{\leftrightarrow} b_{i}\biggl) \bigg],
  \notag
\end{eqnarray}
where the first sum is over all possible locations of first and
N-th regeneration points, and the second sum is over all possible sets of
increments between $v$ and $w$. If the magnitude of one of these
increments is greater than $2d$, this implies at least one
regeneration point must be outside the tube $T_{d}^{i}$.
Therefore, we can restrict the second sum to $\zeta' \in V(v,w)
\cap [-2d,2d]^{N}$, and the last sum is bounded by
\begin{equation} \label{E:boxrestric}
  \sum_{v,w\in T_{d}^{i}\cap \mathcal{S}^{t_{i}}_{a_{i}+3e_i,b_{i}}
  \cap (\mathbb{Z}^{2})^*} \bigg[ \sum_{\zeta' \in V \cap [-2d,2d]^{N}}
  P\biggl(\{a_{i} \overset{\widetilde{h_{t_{i}}}}{\leftrightarrow}
  b_{i} \mbox{ in } T_{d}^{i}\} \cap E\cap F(v,w,\zeta') \ \bigl| \
  a_{i} \overset{\widetilde{h_{t_{i}}}}{\leftrightarrow}
  b_{i}\biggl) \bigg].
\end{equation}
For the remainder of the proof, our sums are over $v,w\in
T_{d}^{i}\cap \mathcal{S}^{t_{i}}_{a_{i}+3e_i,b_{i}} \cap
(\mathbb{Z}^{2})^*$ and $\zeta' \in V(v,w) \cap [-2d,2d]^{N}$. We
can write the last expression as
\begin{eqnarray*}
& &\sum_{v,w}\bigg[\sum_{\zeta'} P\bigl(E \cap F(v,w,\zeta') | \ a_{i}
  \overset{\widetilde{h_{t_{i}}}}{\leftrightarrow} b_{i}\bigl) \
  P\big(a_{i}
  \overset{\widetilde{h_{t_{i}}}}{\leftrightarrow} b_{i} \mbox{ in }
  T_{d}^{i}\ \bigl|  \ E\cap F(v,w,\zeta')\big) \bigg] \\
&\leq& 2 \exp \biggl(\frac{-\kappa|w_{i+1}-w_{i}|}{d^{2}}\biggl) 
  \sum_{v,w}\sum_{\zeta'} P\bigl(E
  \cap F(v,w,\zeta') | \ a_{i}
  \overset{\widetilde{h_{t_{i}}}}{\leftrightarrow} b_{i}\bigl),
\end{eqnarray*}
using~$\eqref{E:intmainlemma}$. Taking the double sum over the
probabilities of disjoint events, in view of \eqref{E:vwsum} and
\eqref{E:boxrestric} we get
\begin{eqnarray*}
  P\biggl(\{a_{i}
  \overset{\widetilde{h_{t_{i}}}}{\longleftrightarrow} b_{i} \mbox{
  in } T_{d}^{i}\} \cap E \ \bigl| \ a_{i}
  \overset{\widetilde{h_{t_{i}}}}{\longleftrightarrow} b_{i}\biggl)
  \leq 2 \exp \biggl(\frac{-\kappa|w_{i+1}-w_{i}|}{d^{2}}\biggl).
\end{eqnarray*}
Combining this with$~\eqref{E:Ec}$ completes the proof.
\end{proof}

\section{Assembling the Segments}
In the last section, we proved that on every long facet of the
$\partial$HSkel$_{s}(\Gamma_{0})$, there is an extra probabilistic
cost for staying in the narrow tube. In this section, we will
bring the pieces in the preceding sections together to deduce
that, there is an extra probabilistic cost of staying in the
annular region $A_{d}$ (with diameter $2d$), throughout the
boundary of the HSkel$_{s}(\Gamma_{0})$. We will show that in light
of the inequality $\eqref{E:skeltohull}$, leaving the annular region
$A_{d}$ implies that MLR$(\Gamma_0) > \theta l^{1/3}(\log l)^{-2/3}$. 
Then
by bounding the number of possible skeletons, we will prove
$\eqref{E: mainthm}$.

\begin{Lemma} There exists $K_{20}=K_{20}(\delta,\gamma)$ such that for
sufficiently large $l$, for all $l$-regular $s$-hull skeletons $\{
w_0,..,w_{m+1} \}$,
\begin{eqnarray}
  \label{E:piecing}\\ \nonumber P\big( 
  \big\{w_{0}\leftrightarrow w_{1}\big\} \circ \dots \circ
  \big\{w_{m} \leftrightarrow w_{m+1}\big\} \mbox{ in }A_{d}\big) \leq
  \exp(-\mathcal{W}_1l-\frac{K_{20}}{\theta ^{2}}l^{1/3}(\log l)^{4/3}).
\end{eqnarray}
\end{Lemma}

\begin{proof}
Using the BK-inequality, we have
\begin{eqnarray} \label{E:BKuse}
  P\big( 
  \big\{w_{0}\leftrightarrow w_{1}\big\} \circ \dots \circ
  \big\{w_{m} \leftrightarrow w_{m+1}\big\} \mbox{ in }A_{d}\big)
  \leq \prod_{i=0}^{m} P(w_{i} \leftrightarrow w_{i+1} \mbox{ in }
  A_{d}).
\end{eqnarray}
This last product can be written as products over
long and short sides separately. We will bound the product over
long sides further. As before for $i\in \mathcal{L}$, let $a_{i} \in
\mathcal{S}^{t_{i}}_{w'_{i},w'_{i}+e_i} \cap T_{d}^{i}$ and
$b_{i}\in \mathcal{S}^{t_{i}}_{w''_i,w''_i-e_i} \cap
T_{d}^{i}$; note there are at most $2d$ choices each for $a_i$ and
$b_i$.  Using$~\eqref{E: Adtotubelemma}$, $\eqref{E:mainlemma}$ and
$l$-regularity we have
\begin{eqnarray} \label{E: psiover2}
& &\prod_{i \in \mathcal{L}} P(w_{i} \leftrightarrow w_{i+1} \mbox{
  in } A_{d}) \\
& &\leq  \ \prod_{i \in \mathcal{L}} \ P( w_{i} \leftrightarrow
  w_{i+1})\biggl[K_{14}\exp{(-K_{15}l)}+4 d^{2}\max_{a_{i},b_{i}}
  P( a_{i} \leftrightarrow b_{i}
  \mbox{ in } T_{d}^{i} \  \bigl| \    a_{i}
  \overset{\widetilde{h_{t_{i}}}}{\leftrightarrow} b_{i})\biggl]
  \notag \\ 
& &\leq \prod_{i \in \mathcal{L}} \ P( w_{i} \leftrightarrow
  w_{i+1})\biggl[K_{14}\exp{(-K_{15}l)}+12 d^{2}
   \exp{\bigg(\frac{-\kappa |w_{i+1}-w_{i}|}{d^{2}}\bigg)}\biggl]
  \notag\\ 
& &\leq \prod_{i \in \mathcal{L}} \ P( w_{i} \leftrightarrow
  w_{i+1})\biggl[13 d^{2}
   \exp{\bigg(\frac{-\kappa |w_{i+1}-w_{i}|}{d^{2}}\bigg)}\biggl]. \notag
\end{eqnarray}
Now we place a condition on the as-yet-unspecified constant $\theta $;
recall that
\[
  s=\left( \frac{\theta \sqrt{\pi}}{2K_{7}} \right)^{1/2}
  l^{2/3} (\log l)^{-1/3}, \qquad d=2 \theta  l^{1/3} (\log l)^{-2/3}.
\]
For some $\beta=\beta(K_{5},K_{7},\kappa)$, provided $\theta $ is
sufficiently small we have
\begin{eqnarray*}
  13 d^{2}
   \exp{\bigg(\frac{-\kappa |w_{i+1}-w_{i}|}{2d^{2}}\bigg)} \leq 20
  \theta ^{2} l^{2/3} (\log l)^{-4/3} \exp (-\beta \theta ^{-3/2}
  \log l)  \leq 1,
\end{eqnarray*}
so that
\begin{eqnarray*}
  13 d^{2} \exp{\bigg(\frac{-\kappa |w_{i+1}-w_{i}|}{d^{2}}\bigg)} \leq
  \exp{\bigg(\frac{-\kappa |w_{i+1}-w_{i}|}{2d^{2}}\bigg)}.
\end{eqnarray*}
Therefore, using Lemma \ref{L:Lsum}, the right side of $~\eqref{E:
psiover2}$ is bounded by
\begin{eqnarray*}
& &\exp{ \biggl(-\frac{\kappa}{2d^{2}}\  \cdot
  \sum_{i \in \mathcal{L}}|w_{i+1}-w_{i}|\biggl)} \cdot \prod_{i \in
  \mathcal{L}} P(w_{i} \leftrightarrow w_{i+1}) \\
& &\qquad\leq \exp \biggl( -\frac{\kappa\sqrt{\pi}}{8 \theta ^{2}\sqrt{2}}\
  l^{1/3} (\log l)^{4/3}\biggl) \cdot \prod_{i \in \mathcal{L}}
  P(w_{i} \leftrightarrow w_{i+1}).
\end{eqnarray*}
Using (\ref{E: psiover2}) and $\eqref{E:coniq}$ we then obtain
\begin{eqnarray} \label{E:prodbound}
 \prod_{i=0}^{m} P(w_{i} \leftrightarrow w_{i+1} \mbox{ in }
  A_{d}) &\leq& \exp \biggl( -\frac{\kappa\sqrt{\pi}}{8
  \theta ^{2}\sqrt{2}}\  l^{1/3} (\log l)^{4/3}\biggl) \prod_{i=0}^{m}
  P(w_{i}
  \leftrightarrow w_{i+1}) \notag \\
&\leq& \exp \biggl( -\frac{\kappa\sqrt{\pi}}{8 \theta ^{2}\sqrt{2}}\
  l^{1/3} (\log l)^{4/3}-\sum_{i=0}^{m} \tau(w_{i+1}-w_{i})\biggl).
\end{eqnarray}
By $l$-regularity there exists a dual circuit $\gamma_0$ with 
$|\Int(\gamma_0)| \geq l^2$, diam$(\gamma_0) \leq 8 \sqrt{2}\ l$
and HSkel$_s(\gamma_0) = \{w_0,..,w_{m+1}\}$.  The first condition
implies diam$(\gamma_0) \geq l$, and by definition of $\mathcal{W}_1$ we have
$\mathcal{W}(\partial$Co$(\gamma_0) \geq \mathcal{W}_1 l$.  Therefore by 
$\eqref{E:cg4}$, for some $K_{21}$,
\[
  \sum_{i=0}^{m} \tau(w_{i+1}-w_{i}) \geq \mathcal{W}_1 l -
  K_{21}\theta l^{1/3}(\log l)^{-2/3}.
\]
The lemma now follows from this together with (\ref{E:BKuse}) and
(\ref{E:prodbound}).
\end{proof}

\begin{proof}[Proof of Theorem \ref{E: mainthm}:] 
By (\ref{E:Eucldiam}),
\begin{eqnarray}
& & P\big(\ \MLR(\Gamma_{0}) \leq \theta  l^{1/3} (\log
  l)^{-2/3} \ \big|
  \  |\Int(\Gamma_{0})| \geq l^{2}  \ \big)  \notag \\
& & \qquad\leq P\big(\ \{\MLR(\Gamma_{0}) \leq \theta  l^{1/3} (\log
  l)^{-2/3}\}
  \cap
  \{\diam(\Gamma_{0})\leq 8 \sqrt{2}l\}
  \ \big|
  \  |\Int(\Gamma_{0})| \geq l^{2}  \ \big) \label{E:lgdiamout}\\
& & \qquad \qquad + \exp(-\mathcal{W}_1l/2). \notag
\end{eqnarray}
This means we need only consider $l$-regular skeletons
$\{w_{0},\cdots,w_{m+1}\}$ for $\Gamma_0$.  When $|\Int(\Gamma_0)|
\geq l^2$ we have diam$(\Gamma_0) \geq l$ and therefore
\[
  \frac{K_{7} s^2}{\diam(\Gamma_0)} < d.
\]
Presuming HSkel$_s(\Gamma_0) = \{w_{0},\cdots,w_{m+1}\}$, this
implies
$\partial$Co$(\Gamma_0) \subset A_d$.  This means that in order to have
$\MLR(\Gamma_{0}) \leq \theta  l^{1/3} (\log l)^{-2/3}$, $\Gamma_0$
must be entirely inside $A_d$.  Thus
\begin{eqnarray}
& & P\big(\ \{\MLR(\Gamma_{0}) \leq \theta  l^{1/3} (\log
  l)^{-2/3}\}
  \cap
  \{\diam(\Gamma_{0})\leq 8 \sqrt{2}l\} \cap
  \{|\Int(\Gamma_{0})|
  \geq l^{2}\}   \big) \notag \\
& &\notag \quad \leq \sum_{\{w_{0},\cdots,w_{m+1}\}} P\bigg(\
  \big\{\mbox{HSkel}_{s}(\Gamma_{0})=\{w_{0},\cdots, w_{m}\}\big\}\\   
& & \qquad \quad \qquad \qquad \quad \cap
  \bigg\{\  \big\{w_{0}\leftrightarrow w_{1}\big\} \circ \dots \circ
  \big\{w_{m} \leftrightarrow w_{m+1}\big\} \mbox{ in }
   A_{d}(w_0,..,w_{m+1}\  \big\}
  \bigg\} \bigg), \label{E:sumbound})
\end{eqnarray}
where the sum is over all $l$-regular skeletons.  By $\eqref{E:cg1}$ the
number of $l$-regular skeletons is at most
\[
  (K_{22}l^2)^{K_{23}\theta ^{-1/2}l^{1/3}(\log l)^{1/3}}
  \leq \exp(K_{24}\theta ^{-1/2}l^{1/3}(\log l)^{4/3}),
\]
for some $K_{22},K_{23},K_{24}$.
This together
with $\eqref{E:piecing}$ and (\ref{E:sumbound}) gives
\begin{eqnarray*}
& & P\big(\ \{\MLR(\Gamma_{0}) \leq \theta  l^{1/3} (\log
  l)^{-2/3}\}
  \cap
  \{\diam(\Gamma_{0})\leq 8 \sqrt{2}l\} \cap
  \{|\Int(\Gamma_{0})|
  \geq l^{2}\}   \big)\\
& &\qquad \leq \exp\bigg(-\mathcal{W}_1l- \bigg( \frac{K_{20}}{\theta ^{2}} \
  - \frac{K_{24}}{\theta ^{1/2}} \bigg) l^{1/3}
  (\log l)^{4/3} \bigg) 
\end{eqnarray*}
which with (\ref{E:droplowerbd}) and (\ref{E:lgdiamout}) yields
\begin{align*}
  P\big(\ &\MLR(\Gamma_{0}) \leq \theta  l^{1/3} (\log
    l)^{-2/3} \ \big|
    \  |\Int(\Gamma_{0})| \geq l^{2}  \ \big)\\
  &\leq \exp\bigg(-\bigg( \frac{K_{20}}{\theta ^{2}} \
   - \frac{K_{24}}{\theta ^{1/2}} \bigg) l^{1/3}
    (\log l)^{4/3} + K_{11} l^{1/3} (\log l)^{2/3} \bigg)
    + \exp\bigg(-\frac{\mathcal{W}_1 l}{2}\bigg). \notag
\end{align*}
For $\theta >0$, sufficiently small,
the last bound tends to $0$, as $l \to \infty$. 
\end{proof}


\begin{thebibliography}{99}

\bibitem[Al1]{Al1} Alexander, K.S., \textit{Approximation of
subadditive functions and rates of convergence in limiting shape
results}, Ann. Probab. \textbf{25} (1997), 30--55.

\bibitem[Al2]{Al2} Alexander, K.S., \textit{Power-law corrections to
exponential decay of connectivities and correlations in lattice
models}, Ann. Probab. \textbf{29} (2001), 92--122.

\bibitem[Al3]{Al3} Alexander, K.S., \textit{Cube-root boundary
fluctuations for droplets in random cluster models}, Commun. Math
Phys. \textbf{224} (2001), 733--781.

\bibitem[ACC]{ACC} Alexander K.S., Chayes, J.T., and Chayes L.,
\textit{The Wulff construction and asymptotics of the finite cluster
distribution for two dimensional Bernoulli percolation}, Commun.
Math. Phys. \textbf{131} (1990), 1-50.

\bibitem[BDJ]{BDJ} Baik, J., Deift, P. and Johansson K., \textit{On
the distribution of the length of the longest increasing subsequence
of random permutations}, J. Amer. Math. Soc. \textbf{12} (1999),
1119--1178.

\bibitem[CI]{CI} Campanino, M. and Ioffe D., \textit{Ornstein-Zernike
Theory for the Bernoulli bond percolation on $\mathbb{Z}^{d}$},
Ann. Probab. \textbf{30} (2002), 652--682.

\bibitem[DH]{DH} Dobrushin, R.L. and Hryniv, O., \emph{Fluctuations 
of the phase boundary
in the $2D$ Ising ferromagnet}, Commun. Math. Phys. \textbf{189} (1997),
395--445.

\bibitem[DKS]{DKS} Dobrushin, R.L., Koteck\'y, R. and Shlosman, S., 
\emph{Wulff construction. A global shape from local
interaction},  Translations of Mathematical Monographs, 104, American 
Mathematical Society, Providence (1992).

\bibitem[ES]{ES} Edwards, R.G. and Sokal, A.D., 
\emph{Generalization of the
Fortuin-Kasteleyn-Swendsen-Wang representation and Monte Carlo algorithm},
Phys. Rev. D \textbf{38} (1988), 2009-2012.

\bibitem[FK]{FK} Fortuin, C.M. and Kasteleyn, P.W., \emph{On
the random cluster model. I. Introduction and relation to other models},
Physica \textbf{57} (1972), 536-564.

\bibitem[FKG]{FKG} Fortuin, C.M., Kasteleyn, P.W. and Ginibre, J.,
\textit{Correlation inequalities on some partially ordered sets},
Commun. Math. Phys. \textbf{22} (1971), 89--103.

\bibitem[Ha]{Ha} Harris T.E., \textit{A lower bound for the critical
probability in a certain percolation process}, Proc.
Camb. Phil. Soc. \textbf{56} (1960), 13--20.

\bibitem[Ho]{Ho} Hoeffding, W., \emph{Probability inequalities for sums
of bounded random variables}, J. Amer. Statist. Assoc. \textbf{58}
(1953), 13--30.

\bibitem[Hr]{Hr} Hryniv, O., \emph{On local behaviour of the phase 
separation line in the $2D$ Ising model}, Probab. Theory Rel. Fields
\textbf{110} (1998), 91-107.

\bibitem[IS]{IS} Ioffe, D. and Schonmann, R.H.,
\textit{Dobrushin-Kotecky-Shlosman theorem up to the critical
temperature}, Commun. Math Phys. \textbf{199} (1998), 91--107.

\bibitem[Jo]{Jo} Johansson, K., \emph{Discrete orthogonal polynomial
ensembles and the Plancherel measure}, Ann. Math. (2) \textbf{153}
(2001), 259--296.

\bibitem[KPZ]{KPZ} Kardar, M., Parisi, G. and Zhang, Y.-C.,
\emph{Dynamic scaling of growing interfaces}, Phys. Rev. Lett.
\textbf{56} (1986), 889--892.

\bibitem[Ke]{Ke} Kesten, H., \textit{The critical probability of bond
percolation on the square lattice equals $\frac{1}{2}$},
Commun. Math. Phys. \textbf{74} (1980), 41--59.

\bibitem[KS]{KS} Krug, J. and Spohn, H., \emph{Kinetic roughening of 
growing interfaces}, in \emph{Solids Far from Equilibrium:  
Growth, Morphology
and Defects} (C. Godr\`eche, ed.) 479--582, Cambridge University Press,
Cambridge (1991).

\bibitem[LNP]{LNP} Licea, C., Newman, C.M. and Piza, M.S.T., 
\emph{Superdiffusivity in
first-passage percolation}, Probab. Theory Rel. Fields 
\textbf{106} (1996), 
559--591.

\bibitem[Me]{Me} Menshikov, M.V., \emph{Coindidence of critical points
in percolation problems}, Soviet Math. Dokl. \textbf{33} (1986),
856--859.

\bibitem[MS1]{MS1} Minlos, R.A. and Sinai Ya.G.,\textit{The
phenomenon of ``phase separation'' at low temperatures in some lattice
models of a gas. I.}, Mat. Sb. \textbf{73} (1967), 375--448.  [English
transl., Math. USSR-Sb. \textbf{2} (1967), 335--395.]

\bibitem[MS2]{MS2} Minlos, R.A. and Sinai Ya.G., \textit{The
phenomenon of "phase separation" at low temperatures in some lattice
models of a gas. II.}, Tr. Moskov. Mat. Obshch. \textbf{19} (1968),
113--178. [English transl., Trans. Moscow Math Soc. \textbf{19} (1968),
121--196.]

\bibitem[NP]{NP} Newman, C.M. and Piza, M.S.T., 
\emph{Divergence of shape fluctuations
in two dimensions}, Ann. Probab. \textbf{23} (1995), 977--1005.

\bibitem[Pi]{Pi} Piza, M.S.T., \emph{Directed polymers in a random
environment:  Some results on fluctuations}, J. Statist. Phys.
\textbf{89} (1997), 581--603.

\bibitem[Ta1]{Ta1}  Taylor, J.E., \emph{Existence and structure of
solutions to a class of nonelliptic variational problems}, Symp. Math.
\textbf{14} (1974), 499-508.

\bibitem[Ta2]{Ta2} Taylor, J.E., \textit{Unique structure of
solutions to a class of nonelliptic variational problems}, Proc.
Sympos. Pure Math. \textbf{27} (1975), 419--427.

\bibitem[Uz]{Uz} Uzun, H.B., \textit{On maximum local roughness of
random droplets in two dimensions}, Ph.D. dissertation, Univ. of
Southern California (2001).

\bibitem[vdBK]{vdBK} van den Berg, J. and Kesten, H.,
\textit{Inequalities with applications to percolation and
reliability}, Journal of Applied Probability \textbf{22} (1985),
556--569.

\bibitem[Wu]{Wu} Wulff, G., \textit{Zur Frage der Geschwingkeit des
Wachstums und der Aufl\"osung der
Krystallflachen}, Z. Kryst. \textbf{34} (1901), 449-530.

\end{thebibliography}
\end{document}